\documentclass[11pt]{article}
\usepackage{graphicx}
\usepackage[tight]{subfigure}
\usepackage{amssymb}

\newcommand{\comment}[1]{}

\def\squarebox#1{\hbox to #1{\hfill\vbox to #1{\vfill}}}
\def\qed{\hspace*{\fill}
        \vbox{\hrule\hbox{\vrule\squarebox{.667em}\vrule}\hrule}\smallskip}
\newenvironment{proof}{\begin{trivlist}
  \item[\hspace{\labelsep}{\em\noindent Proof.~}]
  }{\qed\end{trivlist}}

\newtheorem{lemma}{Lemma}[section]
\newtheorem{theorem}[lemma]{Theorem}
\newtheorem{corollary}[lemma]{Corollary}

\newtheorem{claim}[lemma]{Claim}
\newtheorem{definition}[lemma]{Definition}
\newtheorem{observation}[lemma]{Observation}

\def\squareforqed{\hbox{\rlap{$\sqcap$}$\sqcup$}}
\def\qed{\ifmmode\squareforqed\else{\unskip\nobreak\hfil
\penalty50\hskip1em\null\nobreak\hfil\squareforqed
\parfillskip=0pt\finalhyphendemerits=0\endgraf}\fi}

\newcommand{\nats}{\mbox{$\mathbb N$}}

\newcommand{\ch}{\mbox{\rm ch}}
\newcommand{\ind}{\mbox{\rm ind}}

\newcommand{\diam}{\mbox{\rm diam}}

\newcommand{\pr}{\mbox{\rm pr}}
\newcommand{\redscale}{0.35}   
\newcommand{\Greedy}{\mbox{\sf Greedy}}

\newcommand{\defn}[1]{#1}
\newcommand{\mathdefn}[1]{$#1$}

\newlength{\tablength}
\newlength{\spacelength}
\newcommand{\tabstar}{\hspace*{\tablength}}
\newcommand{\spacestar}{\hspace*{\spacelength}}
\def\obeytabs{\catcode`\^^I=\active}
{\obeytabs\global\let^^I=\tabstar}
{\obeyspaces\global\let =\spacestar}
\newenvironment{display}{\begingroup\obeylines\obeyspaces\obeytabs}{\endgroup}
\newenvironment{prog}{\begin{display}\parskip0pt\sf}{\end{display}}

\title{On Colorings of Squares of Outerplanar Graphs\thanks{An earlier
version of this current paper appeared in SODA 2004~\cite{GeirMagnus-SODA2}.
}}

\author{
{\sl Geir Agnarsson}
\thanks{Department of Mathematical Sciences,
George Mason University, 
MS 3F2, 
4400 University Drive, 
Fairfax, VA 22030, 
{\tt geir@math.gmu.edu}} 
\and 
{\sl Magn\'{u}s M.~Halld\'{o}rsson}
\thanks{
Department of Computer Science, University of Iceland, Reykjav\'{\i}k, Iceland.
{\tt mmh@hi.is}}}

\begin{document}

\date{}

\maketitle

\begin{abstract}
\noindent
We study vertex colorings of the square $G^2$ of 
an outerplanar graph $G$. We find the optimal bound 
of the inductiveness, chromatic number and the clique
number of $G^2$ as a function of the maximum degree
$\Delta$ of $G$ for all $\Delta\in \nats$. 
As a bonus, we obtain the optimal bound of the 
choosability (or the list-chromatic number) of $G^2$ 
when $\Delta \geq 7$. In the case of chordal 
outerplanar graphs, we classify exactly which graphs have 
parameters exceeding the absolute minimum.

\vspace{3 mm}


\vspace{2 mm}

\end{abstract}


\section{Introduction}
\label{sec:intro}

The square of a graph $G$ is the graph $G^2$ on the same vertex set
with edges between pair of vertices of distance one or two in $G$.
Coloring squares of graphs has been studied, e.g., in relation to frequency
allocation. This models the case when nodes represent both senders and
receivers, and two senders with a common neighbor will interfere if
using the same frequency. 

The problem of coloring squares of graphs has particularly seen much
attention on planar graphs.
A conjecture of Wegner~\cite{Wegner} dating from 1977 (see
\cite{JensenToft}), states that the square of every planar graph $G$ of
maximum degree $\Delta \ge 8$ has a chromatic number which does not
exceed $3\Delta/2 + 1$.  The conjecture matches the maximum clique
number of these graphs.  Currently the best upper bound known is
$1.66\Delta + 78$ by Molloy and Salavatipour~\cite{MS01}.

An earlier paper of the current authors~\cite{GeirMagnus} gave a bound
of $\lceil 1.8\Delta \rceil$ for the chromatic number of squares of
planar graph with large maximum degree $\Delta \ge 749$. This is based
on bounding the {\em inductiveness} of the graph, which is the maximum
over all subgraphs $H$ of the minimum degree of $H$. It was also shown
there that this was the best possible bound on the inductiveness. Borodin et
al~\cite{BBGH01} showed that this bound holds for all $\Delta \ge 48$.
Inductiveness has the additional advantage of also bounding the
list-chromatic number. 

Inductiveness leads to a natural greedy algorithm (henceforth called \Greedy):
Select vertex $u\in V(G)$ of minimum degree, sometimes called 
a {\em simplicial} vertex of $G$, recursively color $G \setminus u$, 
and finally color $u$ with the smallest available color.
Alternatively, $k$-inductiveness leads to an {\em inductive ordering}
$u_1, u_2, \ldots, u_n$ of the vertices such that any vertex $u_i$ has at most $k$
neighbors among $\{u_{i+1}, \ldots, u_n\}$. Then, if we color the
vertices {\em first-fit} in the
reverse order $u_n, u_{n-1}, \ldots, u_1$ (i.e.~assigning each vertex
the smallest color not used among its previously colored neighbors),
the number of colors used is at most $k+1$. Implemented efficiently,
the algorithm runs in time linear in the size of the graph~\cite{mmh:greed}.
The algorithm has also the special advantage
that it requires only the square graph
$G^2$ and does not require information about the underlying graph $G$.

The purpose of this article is to further contribute to the study of
various vertex colorings of squares of planar graphs, by examining an
important subclass of them, the class of outerplanar graphs. 
Observe that the neighborhood of a vertex with $\Delta$ neighbors induces
a clique in the square graph. Thus, the chromatic number, and in fact
the clique number, of any graph of maximum degree $\Delta$ is
necessarily a function of $\Delta$ and always at least $\Delta+1$. 


\paragraph{Our results.}
We derive tight bounds on chromatic number, as well as the
inductiveness and the clique number of the square of an outerplanar graph
$G$ as a function of the maximum degree $\Delta$ of $G$.
One of the main results, given in Section~\ref{sec:ind}, is that when $\Delta \ge 7$, the
inductiveness of $G^2$ is exactly $\Delta$. It follows that the clique
and chromatic numbers are exactly $\Delta+1$ and that {\Greedy} yields
an optimal coloring. As a bonus we obtain in this case that the 
{\em choosability} (see Definition~\ref{def:choose}) is 
the optimal $\Delta+1$. We can then treat the low-degree cases
separately to derive a linear-time algorithm independent of $\Delta$.
We examine in detail the low-degree cases, $\Delta < 7$, and
derive best possible upper bounds on the maximum clique and chromatic numbers,
as well as inductiveness of squares of outerplanar graphs.
These bounds are illustrated in Table~\ref{tab:Delta}.
We treat the special case of chordal outerplanar graphs separately,
and further classify all chordal outerplanar graphs $G$ for
which the inductiveness of $G^2$ exceeds $\Delta$ or the clique or
chromatic number of $G^2$ exceed $\Delta+1$.
\begin{table}[htbp]
\begin{center}
\begin{tabular}{|l||l|l|l||l|l|l||}
\label{tab:Delta}
& \multicolumn{3}{c||}{Chordal} & \multicolumn{3}{c||}{General} \\
\hline 
$\Delta$ & $\omega$    & $\ind$    & $\chi$
   & $\omega$    & $\ind$    & $\chi$ \\ \hline
2        & $\Delta + 1$   & $\Delta$     & $\Delta + 1$
         & $\Delta + 3$   & $\Delta+2$   & $\Delta + 3$ \\
3        & $\Delta + 1$   & $\Delta$     & $\Delta + 1$
         & $\Delta + 2$   & $\Delta+1$   & $\Delta + 2$ \\
4        & $\Delta + 2$   & $\Delta +1$  & $\Delta + 2$
         & $\Delta + 2$   & $\Delta +2$  & $\Delta + 2$ \\
5        & $\Delta + 1$   & $\Delta +1$  & $\Delta + 1$
         & $\Delta + 1$   & $\Delta +1$  & $\Delta + 2$\hspace*{-1ex} \\
6        & $\Delta + 1$   & $\Delta +1$  & $\Delta + 1$
         & $\Delta + 1$   & $\Delta +1$  & $\Delta + 1$ \\
7+\hspace*{-1ex}   & $\Delta + 1$   & $\Delta$     & $\Delta + 1$
         & $\Delta + 1$   & $\Delta$     & $\Delta + 1$ \\
\hline 
\end{tabular}
\caption{Optimal upper bounds for the clique number, inductiveness, 
and chromatic number of the square of a chordal / non-chordal outerplanar graph $G$.}
\end{center}
\end{table}

\smallskip


\paragraph{Related results.}
It is straightforward to show that the inductiveness of a square graph
of an outerplanar graph of degree $\Delta$ is at most $2\Delta$ (see
\cite{GeirMagnus}), and this is attained by an inductive ordering of $G$. 
Calamoneri and Petreschi~\cite{CalaPet} gave a linear time
algorithm to distance-2 color outerplanar graphs, as well as for
related problems.
They showed that it uses an optimal $\Delta+1$ colors whenever $\Delta
\ge 7$, and at most $\Delta+2$ colors for $\Delta \ge 4$.
In comparison, we give tight upper and lower bounds for all values of
$\Delta$, give a thorough treatment of the subclass of chordal graphs,
and analyze a generic parameter, inductiveness, that gives as a bonus similar bounds
for the list chromatic number.

Zhou, Kanari and Nishizeki~\cite{ZKN00} gave a polynomial time
algorithm to find an optimal coloring of any power of a partial
$k$-tree $G$, given $G$. Since outerplanar graphs are partial 2-trees,
this solves the coloring problem we consider. For squares of
outerplanar graphs, their algorithm has complexity
$O(n (\Delta+1)^{2^{37}} + n^3)$, 
which is impractical for any values of $\Delta$ and $n$.
When $\Delta$ is constant, one can use the observation of 
Krumke, Marathe and Ravi~\cite{KMR:dialm} that squares of outerplanar
graphs have treewidth at most $k \le 3\Delta-1$.
Thus, one can use efficient ($2^{O(k)}n$ time) algorithms for coloring
partial $k$-trees, obtaining a linear-time algorithm when $\Delta$ is constant.


\paragraph{Organization.}
The rest of the paper is organized as follows:
In Section~\ref{sec:prelim} we introduce our notation and definitions, and
show how the problems regarding the clique number and chromatic number
reduce to the case of biconnected outerplanar graphs.
Inductiveness is treated in Section~\ref{sec:ind}. We then
treat the chordal case in Section~\ref{sec:chordal}. Many examples
here show that the lower bounds derived in other sections 
(i.e.~Sections~\ref{sec:ind} and~\ref{sec:clique}) are 
optimal. The clique number is derived in Section~\ref{sec:clique}.
The last Section~\ref{sec:col} derives optimal bounds on chromatic number 
in each of the smaller cases of $\Delta \in\{2,3,4,5,6\}$. The main
result there is the optimal bound for the
chromatic number of $G^2$ in the hardest case when $\Delta = 6$.


\section{Definitions}
\label{sec:prelim}

In this section we give some basic definitions and prove results that will 
be used later for our results in the following sections.


\paragraph{Graph notation.}
The set $\{1,2,3,\ldots \}$ of natural numbers will be denoted
by $\nats$. Unless otherwise stated, a graph $G$ will always be
a simple graph $G = (V,E)$ where $V = V(G)$ is the set of vertices
or nodes, and $E = E(G)$ the set of edges of $G$. The edge between the vertices
$u$ and $v$ will be denoted by $uv$ (here $uv$ and $vu$ will mean the
same undirected edge). 
By {\em coloring} we will always mean vertex coloring.
We denote by \mathdefn{\chi(G)} the chromatic number of $G$ and
by \mathdefn{\omega(G)} the clique number of $G$.
The degree of a vertex $u$ in graph $G$ is denoted by  $d_G(u)$.
We let \mathdefn{\delta(G)} and \mathdefn{\Delta(G)} denote 
the minimum and maximum 
degree of a vertex in $G$ respectively. 
When there is no danger of ambiguity, we simply write $\Delta$ instead of $\Delta(G)$. 
We denote by \mathdefn{N_G(u)} the open neighborhood of $u$ in $G$, that is the set
of all neighbors of $u$ in $G$, and by $N_G[u]$ the closed neighborhood 
of $u$ in $G$, that additionally includes $u$.

The square graph \mathdefn{G^2} of a graph $G$ is a graph on the same vertex set
as $G$ in which additionally to the edges of $G$, 
every two vertices with a common neighbor in $G$ are
also connected with an edge. Clearly this is the same as the graph
on $V(G)$ in which each pair of vertices of distance 2 or less in $G$
are connected by an edge in $G^2$. 

By a \defn{{\em $k$-vertex}} we will mean a vertex of degree at most 2 in $G$ 
and distance-2 degree at most $k$. 
`

\paragraph{Tree terminology.} 
The {\em diameter} of $T$ is the number of edges in the longest simple path in 
$T$ and will be denoted by \mathdefn{\diam(T)}.
For a tree $T$ with $\diam(T)\geq 1$ we can form the {\em pruned} tree 
\mathdefn{\pr(T)} by removing all the leaves of $T$.
A {\em center} of $T$ is a vertex of distance
at most $\lceil \diam(T)/2\rceil$ from all other vertices of $T$.
A center of $T$ is either unique or one of two unique adjacent vertices.
When $T$ is rooted at $r\in V(T)$, the \defn{{\em $k$-th ancestor}}, if it exists,
of a vertex $u$ is the vertex on the unique path 
from $u$ to $r$ of distance $k$ from $u$.
An {\em ancestor} of $u$ is a $k$-th ancestor of $u$ for some $k\geq 0$. 
Note that $u$ is viewed as an ancestor of itself. 
The \defn{{\em parent}} (\defn{{\em grandparent}}) of a vertex is then the
$1$-st ($2$-nd) ancestor of the vertex. The \defn{{\em sibling}} of a vertex is
another child of its parent, and a \defn{{\em cousin}} is child of a sibling
of its parent. 
The \emph{height} of a rooted tree is 
the length of the longest path from the root to a leaf. 
The \emph{height} of a vertex $u$ in a rooted tree $T$ is the height of the 
rooted subtree of $T$ induced by all vertices with $u$ as an ancestor. 
A tree is said to be {\em full} if it contains no degree-two vertices.

Note that in a rooted tree $T$, vertices of height zero are the leaves 
(provided that the root is not a leaf). Vertices of height one are 
the parents of leaves, that is, the leaves of the pruned
tree $\pr(T)$ and so on. In general, for $k\geq 0$ let $\pr^k(T)$ be 
given recursively by $\pr^0(T) = T$ and $\pr^k(T) = \pr(\pr^{k-1}(T))$.
Clearly $V(T)\supset V(\pr(T))\supset\cdots\supset V(\pr^k(T))\supset\cdots$
is a strict inclusion. With this in mind we have an alternative 
``root-free'' description of the height of vertices in a tree. 
\begin{observation}
\label{obs:pruned-ss}
Let $T$ be a tree and $0\leq k\leq \lfloor \diam(T)/2\rfloor$.
The vertices of height $k$ in $T$ are precisely the leaves
of $\pr^k(T)$.
\end{observation}


\paragraph{Inductiveness.}
The {\em inductiveness} or the {\em degeneracy} of a graph $G$, denoted
by \mathdefn{\ind(G)}, is the natural number defined by
\[ 
\ind(G) = \max_{H\subseteq G}\left\{ \delta(H)\right\}, 
\]
where $H$ runs through all the induced subgraphs of $G$. 
If $k\geq \ind(G)$ then we say that $G$ is \defn{{\em $k$-inductive}}.

In a graph $G$ of maximum degree at most $\Delta$, note that for any $u\in V(G)$, 
the vertex set $N_G[u]$ will induce a 
clique in $G^2$, and hence $\omega(G^2), \chi(G^2)\geq \Delta + 1$.
Since $\ind(G^2) + 1 \geq \chi(G^2)$, the upper bound of $\ind(G)$ is
necessarily an increasing function $f(\Delta)$ of $\Delta\in \nats$.
In general, the inductiveness of a graph $G$ yields an ordering $\{ u_1,u_2,\ldots
,u_n\}$ of the vertex set $V(G)$ of $G$, such that each vertex $u_i$
has at most $\ind(G)$ neighbors among the previous vertices
$u_1,\ldots, u_{i-1}$ that is to say $|N_G(u_i)\cap \{u_1,\ldots,
u_{i-1}\}|\leq \ind(G)$.  This gives us an efficient way to color
every graph $G$ by at most $\ind(G) +1$ colors in a greedy fashion. 


\paragraph{Biconnectivity.}
The \emph{blocks} of a graph $G$ are the maximal biconnected subgraphs of $G$.
A \emph{cutvertex} is a vertex shared by two or more blocks.
A \emph{leaf block} is a block with only one cutvertex (or none, if the
graph is already biconnected).

We show here that we can assume, without loss of generality, that $G$
is biconnected when considering the chromatic number or the clique number
of $G^2$: Let $G$ be a graph and $\mathcal{B}$ the set of its biconnected blocks.
In the same way that $\omega(G) = \max_{B\in\mathcal{B}}\{\omega(B)\}$ and  
$\chi(G) = \max_{B\in\mathcal{B}}\{\chi(B)\}$, 
we have the following.
\begin{lemma}
\label{lmm:biconn}
For a graph $G$ with a maximum degree $\Delta$ and
set $\mathcal{B}$ of biconnected blocks we have 
\begin{eqnarray*}
\omega(G^2) & = & \max\{\max_{B\in\mathcal{B}}\{\omega(B^2)\},\Delta+1\}, \\
\chi(G^2)   & = & \max\{\max_{B\in\mathcal{B}}\{\chi(B^2)\},\Delta+1\}.
\end{eqnarray*}
Further, optimal $\chi(B^2)$-colorings of the squares of all the blocks
$B^2$ can be modified to a $\chi(G^2)$-coloring of $G^2$ 
in a greedy fashion.
\end{lemma}
\begin{proof}
First note that a clique of $G^2$ with vertices contained in more
than one block of $G$ must contain the cutvertex of two blocks. Therefore
the clique must be induced by the closed neighborhood of this cutvertex,
and hence of size at most $\Delta +1$. This proves the first
formula for $\omega(G^2)$.

For the chromatic number of $G^2$, we proceed by induction on $b =
|\mathcal{B}|$. 
The case $b=1$ is a tautology, so assume $G$ has $b\geq 2$ 
blocks and that the lemma is true for
$b-1$. Let $B$ be a leaf block and let 
$G'= \cup_{B'\in\mathcal{B}\setminus\{B\}}B'$,
with $w = V(B)\cap V(G')$ as a cutvertex.
If $\Delta'$ is the maximum degree of $G'$, then by induction hypothesis
$\chi(G'^2) = \max\{\max_{
B'\in\mathcal{B}\setminus\{B\}}\{\chi(B'^2)\},\Delta'+1\}$. 
Assume we have a $\chi(G'^2)$-coloring of $G'^2$ and a $\chi(B^2)$-coloring of 
$B^2$, the latter given by a map $c_B : V(B)\rightarrow \{1,\ldots,\chi(B^2)\}$.
Since $w$ is a cutvertex we have a partition $N_G[w] = \{w\}\cup N_B \cup N_{G'}$,
where $N_B = N_G(w)\cap V(B)$ and $N_{G'} = N_G(w)\cap V(G')$. In the given 
coloring $c_B$ all the vertices in $N_B$ have received distinct colors,
since they all have $w$ as a common neighbor in $B$. Since 
$|N_G[w]|\leq \Delta+1$ there is a permutation
$\sigma$ of $\{1,\ldots,\max\{\chi(B^2),\Delta+1\}\}$ such that 
$\sigma\circ \mathbf{i}\circ c_B$ yields a new $\chi(B^2)$-coloring of $B^2$ 
such that all vertices in 
$N_G[w]$ receive distinct colors (here $\mathbf{i}$ is the inclusion map
of $\{1,\ldots,\chi(B^2)\}$ in 
$\{1,\ldots,\max\{\chi(B^2),\Delta+1\}\}$.)
This together with the given 
$\chi(G'^2)$-coloring of $G'^2$ provides a vertex coloring of $G^2$
using at most 
$\max\{\max\{\chi(B^2),\Delta+1\},\chi(G'^2)\}\leq 
\max\{\max_{B\in\mathcal{B}}\{\chi(B^2)\},\Delta+1\}$,
which completes our proof.
\end{proof}
Note that Lemma~\ref{lmm:biconn} provides a way to extend
distance-2 colorings of the blocks of $G$ to a distance-2 coloring of
the whole of $G$. Thus, 
by Lemma~\ref{lmm:biconn} we can assume our graphs are biconnected,
both when considering clique and chromatic numbers of $G^2$.

For the inductiveness of $G^2$, such an extension property
as Lemma~\ref{lmm:biconn}, to express $\ind(G^2)$ directly
in terms of $\Delta$ and the inductiveness of the blocks of $G$,
is not as straightforward although it can be done. 
This is mainly because the simplicial vertex of a biconnected block 
could be a cut-vertex of the graph. 
We will consider this better in Section~\ref{sec:ind}.
 

\paragraph{Duals of outerplanar graphs.}
For our arguments to come we need a few properties about outerplanar
graphs, the first of which is an easy exercise (See~\cite{West}).
\begin{claim}
\label{clm:deg=2}
Every biconnected outerplanar graph has at least two vertices of degree 2.
\end{claim}
To analyze the inductiveness of an outerplanar graph $G$,
it is useful to consider the {\em weak dual} of $G$, denoted
by \mathdefn{T^*(G)} and defined in the following: 
\begin{lemma}
\label{lmm:dual-out}
Let $G$ be an outerplanar graph with an embedding in the plane.
Let $G^*$ be its geometrical dual, and let
$u^*_{\infty} \in V(G^*)$ be the vertex corresponding
to the infinite face of $G$. Then the {\em weak dual} graph 
$T^*(G) = G^* - u^*_{\infty}$ is a forest
which satisfies the following:
\begin{enumerate}
  \item $T^*(G)$ is tree iff $G$ is biconnected.
  \item $T^*(G)$ has maximum degree at most three, if $G$ is chordal.
\end{enumerate}
\end{lemma}
Note that for a biconnected chordal graph $G$,
there is a one-to-one correspondence $u\leftrightarrow u^*$ 
between the degree-2 vertices $u$ of $G$, and the leaves $u^*$ 
of $T^*(G)$.
\begin{proof} (Lemma~\ref{lmm:dual-out})
 This follows easily by Claim~\ref{clm:deg=2} and induction on $n= |V(G)|$.
\end{proof}
Note that any biconnected chordal outerplanar graph on $n$ vertices can be
constructed in the following way: Start with two
vertices $v$ and $w$ and connect them with an 
edge. For $i = 1$ to $n-2$, inductively connect
a vertex $u_i$ to two endvertices of an
edge which bounds the infinite face. Hence,
after the $i$-th step, the vertex $u_i$ is of degree 2.
Simultaneously we construct the weak dual
tree $T^*(G)$ on the vertices $u_1^*,\ldots ,u_{n-2}^*$,
by adding $u_i^*$ as a leaf to the vertex in $\{u_1^*,\ldots,u_{i-1}^*\}$
corresponding to the face containing the two neighbors of $u_i$ after
the $i$-th step. Hence, we have the following.
\begin{observation}
\label{obs:1-1-corr}
For a biconnected chordal outerplanar graph $G$, there
are two vertices $v,w\in V(G)$ such that 
there is a bijection 
$V(G)\setminus \{v,w\}\rightarrow V(T^*(G))$, given
by $u\mapsto u^*$, such that degree-2 vertices of $G$
correspond to leaves of $T^*(G)$. Further, successfully
removing degree-2 vertices from $G$ will result in removing
leaves from $T^*(G)$ in such a way that the mentioned
correspondence will still hold between degree-2 vertices of 
the altered graph $G$ and the leaves of the altered tree.
\end{observation}
By Lemma~\ref{lmm:dual-out}, $T^*(G)$ for a chordal graph $G$
is a tree of maximum degree 3, and hence each of its leaves has 
at most one sibling.

Note, however, that if $G$ is not chordal then the assignment 
$u\mapsto u^*$ is only surjective and not bijective. 
Both in the chordal and non-chordal case we will call the vertex 
$u^*$ the {\em dual vertex} of $u$. For the non-chordal case,
such a construction can be done in a similar fashion, 
except that we inductively add a path of length 
$\geq 2$ instead of length exactly 2 


\paragraph{Faces and dual leafs.}
For an outerplanar plane graph $G$ two faces of $G$ are said to 
be \defn{{\em adjoint}} (shortened as \emph{adj.}) if they share a common vertex.
A \defn{{\em $k$-face}} is a face $f$ with $k$ vertices and $k$ edges. This 
will be denoted by $|f| = k$.

For a bounded face $f$ of $G$ the corresponding dual vertex of  $T^*(G)$
will be denoted by $f^*$. Note for a chordal $G$ and if $f$ has two 
bounding edges bounding the infinite face then $f^* = u^*$,
the dual vertex of $f$ from above. We will, however, 
speak interchangeably of a face $f$ and its corresponding 
dual vertex $f^*$ (or $u^*$ from above in the chordal case) from $T^*(G)$,
when there is no danger of ambiguity, and we will apply standard 
forest/tree vocabulary to faces from
the tree terminology given previously when each component 
from $T^*(G)$ is rooted at a center.
A \defn{{\em sib}} of a face $f$ is a sibling in $T^*(G)$ that
is adjoint to $f$. 

A face $f$ is \emph{$i$-strongly simplicial}, or \defn{{\em $i$-ss}} for short, if
either $f$ is isolated (that is $G$ consists of $f$ alone),
or $f$ is a leaf in $T^*(G)$ satisfying one of the 
following: (i) $i=0$, or (ii) the parent face of $f$ in $T^*(G)$ is 
$(i-1)$-ss in $\pr(T^*(G))$. Thus, e.g.\ all leafs are $0$-ss, while 
those leafs whose siblings have no children are also 1-ss.


\section{Inductiveness}
\label{sec:ind}

In this section we will derive optimal bounds on inductiveness.
The following is the main result of this section.
\begin{theorem}
\label{thm:main-result}
For an outerplanar graph $G$ of maximum degree $\Delta\geq 5$,
we have $\ind(G^2)\leq \Delta + 1$. 
If further, $\Delta \geq 7$, then $\ind(G^2) = \Delta$.
\end{theorem}
To bound the inductiveness, it is sufficient to show that there always
exists a vertex that has both small degree and small distance-2 degree.
Recall that a $k$-vertex is a vertex of degree at most 2 in $G$ 
and distance-2 degree at most $k$. 
\begin{lemma}
\label{lmm:ind-leq-k}
Suppose any outerplanar graph of maximum degree $\Delta$ contains a
$k$-vertex. Then, any outerplanar graph $G$ of maximum degree $\Delta$ satisfies
$\ind(G^2) \le k$.
\end{lemma}
\begin{proof}
We show this by induction on $|V(G)|$. Let $G$ be an outerplanar graph
with a $k$-vertex $u$.
We choose one incident edge $uv$ and 
form the \emph{contraction} $G/uv$; this 
is the simple graph obtained from $G$ by contracting the edge $uv$ into a single 
vertex $v'$ and keeping all edges that were incident on either 
$u$ or $v$ (deleting multiple copies). Formally, $G/uv$ has the vertex set 
$V(G/uv) = (V(G) \setminus \{u,v\}) \cup \{v'\}$ and edge set
$E(G/uv) = E(G[V \setminus \{u,v\}]) \cup 
\{v' w : uw\in E(G) \mbox{ or } vw \in E(G)\}$.
The set of distance-2 neighbors 
of $v'$ in $G/uv$ properly contains the
set of distance-2 neighbors of $v$ in $G$. Hence, a
$k$-inductive ordering of $(G/uv)^2$ also gives a $k$-inductive ordering of
$G^2$ excluding $u$ and where $v$ is replaced with $v'$. 
Further, since $u$ was of degree at most 2, 
the degree of $v'$ is at most that of $v$, and hence the 
maximum degree of $G/uv$ is at most that of $G$.
By induction, there is such a $k$-inductive ordering of $(G/uv)^2$.
By prepending $u$ to that ordering, replacing $v'$ by $v$, we
obtain a $k$-inductive ordering of $G^2$.
\end{proof}
Recall that a leaf block $B$ of $G$ contains just one cutvertex. 
Call a block $B$ {\em simple} if $\diam(T^*(B))\leq 2$,
that is, $T^*(B)$ is empty, a single vertex, a single edge,
or a star on three or more vertices.
\begin{lemma}
\label{lmm:simple-leafblocks} 
Any simple leaf block contains a $(\Delta+1)$-vertex of $G$
if $\Delta \geq 5$ and a $\Delta$-vertex if $\Delta \geq 6$.
\end{lemma}
\begin{proof}
If $B$ is a single edge, then the leaf node is a $\Delta$-vertex.
If $B$ is a 3-cycle, then either of the non-cut-vertices are $\Delta$-vertices.
When $B$ is a $k$-cycle, for $k \ge 4$, then any node on the
cycle that is not adjacent to the cut vertex is a $4$-vertex.

Assume now that $T^*(B)$ is a single edge or a star on $\geq 3$ vertices.
Clearly we have $\Delta(B)\leq 4$. If $B$ contains a degree-2 vertex $v$ 
of distance 2 or more from the cutvertex, then $v$ is a 6-vertex.
So $v$ is a $(\Delta+1)$-vertex if $\Delta\geq 5$ and $\Delta$-vertex
if $\Delta\geq 6$. If all the degree-2 vertices of $B$ are adjacent
to the cutvertex of $G$, then $B$ is a diagonalized $C_k$ where $k\in\{4,5\}$.
In this case both the degree-2 vertices of $B$ are $\Delta$-vertices of 
$G$. Hence we have the lemma.
\end{proof}
By Lemma~\ref{lmm:simple-leafblocks} we will focus on non-simple 
leaf blocks for the rest of this section.

We start our search for a $k$-vertex at a face that has some nice properties.
Recall the definition of a $i$-ss face. Notice that a face is $1$-ss
iff its parent in $T^*(G)$ has no grandchildren, while it is $2$-ss if
it either has no grandparent in $T^*(G)$, or if its grandparent is not the 3-rd
ancestor of another face. Note that if $B$ is not simple, then for any rooting 
of $T^*(B)$ there are always faces with parents and grandparents.
\begin{lemma} 
\label{lmm:good-face} 
Let $B$ be a non-simple leaf block of $G$ and let $i$ be a
non-negative integer. 
Then, $B$ contains an $i$-ss face $f$, its parent $f'$, its
grandparent $f''$ and an edge $e$ on the boundary of $f''$, such that
$e$ separates $f''$ and all its descendants (its children and
grandchildren) from the rest of $G$.
\end{lemma}
\begin{proof} 
If $G = B$ is biconnected let $r$ be any face. Otherwise, if $G\ne B$, 
let $r$ be any face that contains the cutvertex on its boundary and an
edge bounding the infinite face.
Let $T^*(B)$ be rooted at $r$ and let $f$ be a face of maximum 
distance from $r$ in $T^*(B)$. It is not hard to see that
the center(s) of $T^*(B)$ is(are) on the path from $f$ to $r$,
and hence $f$ is an endvertex of a maximum length path of $T^*(B)$.
By definition of $i$-ss, $f$ is therefore $i$-ss for any $i\geq 0$.
Also, $f$ has a parent $f'$ and a grandparent $f''$.

If $f$ has a great-grandparent $f'''$, then we let $e$ be the edge
that separates $f''$ from $f'''$ in $T^*(B)$. This edge separates
$f''$ and all its descendants from the rest of $B$, including $r$
and thus necessarily also from the rest of $G$.

If $f$ has no great-grandparent, then $r = f''$ and we let $e$ be
an edge incident on the cutvertex and the infinite face.
\end{proof}
\begin{claim} 
\label{clm:ab-cd}
Assume that we have faces $f$, $f'$ and $f''$ in a block $B$ of $G$
as promised by Lemma~\ref{lmm:good-face}.
\begin{enumerate}
  \item If $ab$ is the edge separating $f'$ from $f''$, then either
$a$ or $b$ has degree at most 6 in $G$.
  \item If $cd$ is the edge separating $f$ from $f'$, then either
$c$ or $d$ has degree at most 4 in $G$.
\end{enumerate}
\end{claim}
\begin{proof}
Note that the cutvertex is either neither of the endvertices
of the edge $e$ that separates $f''$ and all its descendants 
from the rest of $G$, or one of them.

The first statement is true since $f$ is 2-ss in the block $B$, 
and the second statement is true since $f$ is 1-ss in $B$.
\end{proof}

\paragraph{Reducible configurations.}
A \emph{configuration} is an induced plane subgraph with certain
vertices specially marked as having no neighbors outside the subgraph.
A configuration is {\em $k$-reducible}\label{page:reduce} for an
integer $k$, if there exists a $k$-vertex for it.  When $k$ is
understood, we shall simply speak of a \emph{reducible configuration}.

\begin{figure} 
\centering
\subfigure[$(\Delta+1)$-inductiveness, $\Delta \ge 5$]{
\label{fig:subfig:d1-flow} 
\includegraphics[scale=0.45]{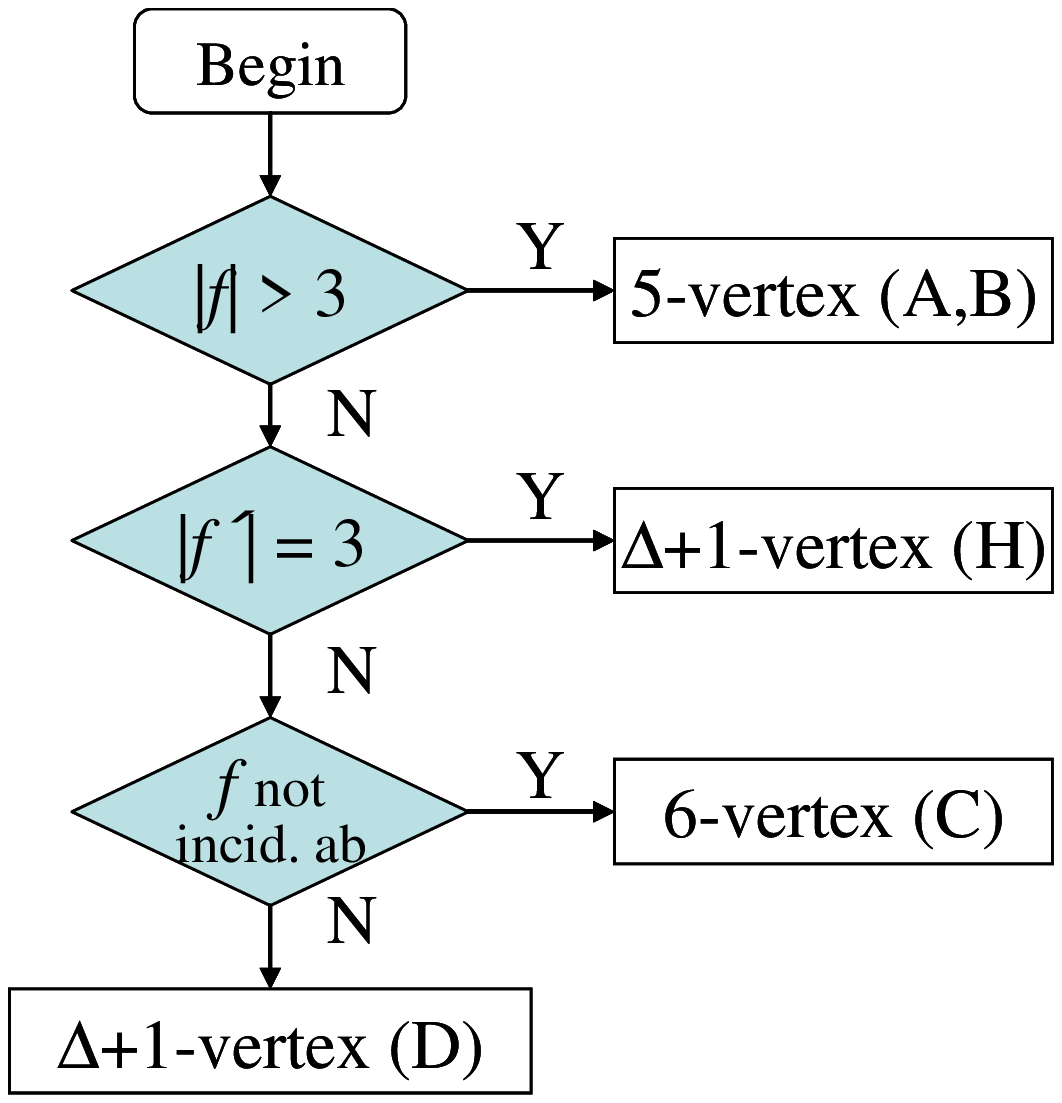}}
\hspace{1in}
\subfigure[$\Delta$-inductiveness, $\Delta \ge 7$]{
\label{fig:subfig:d-flow} 
\includegraphics[scale=0.45]{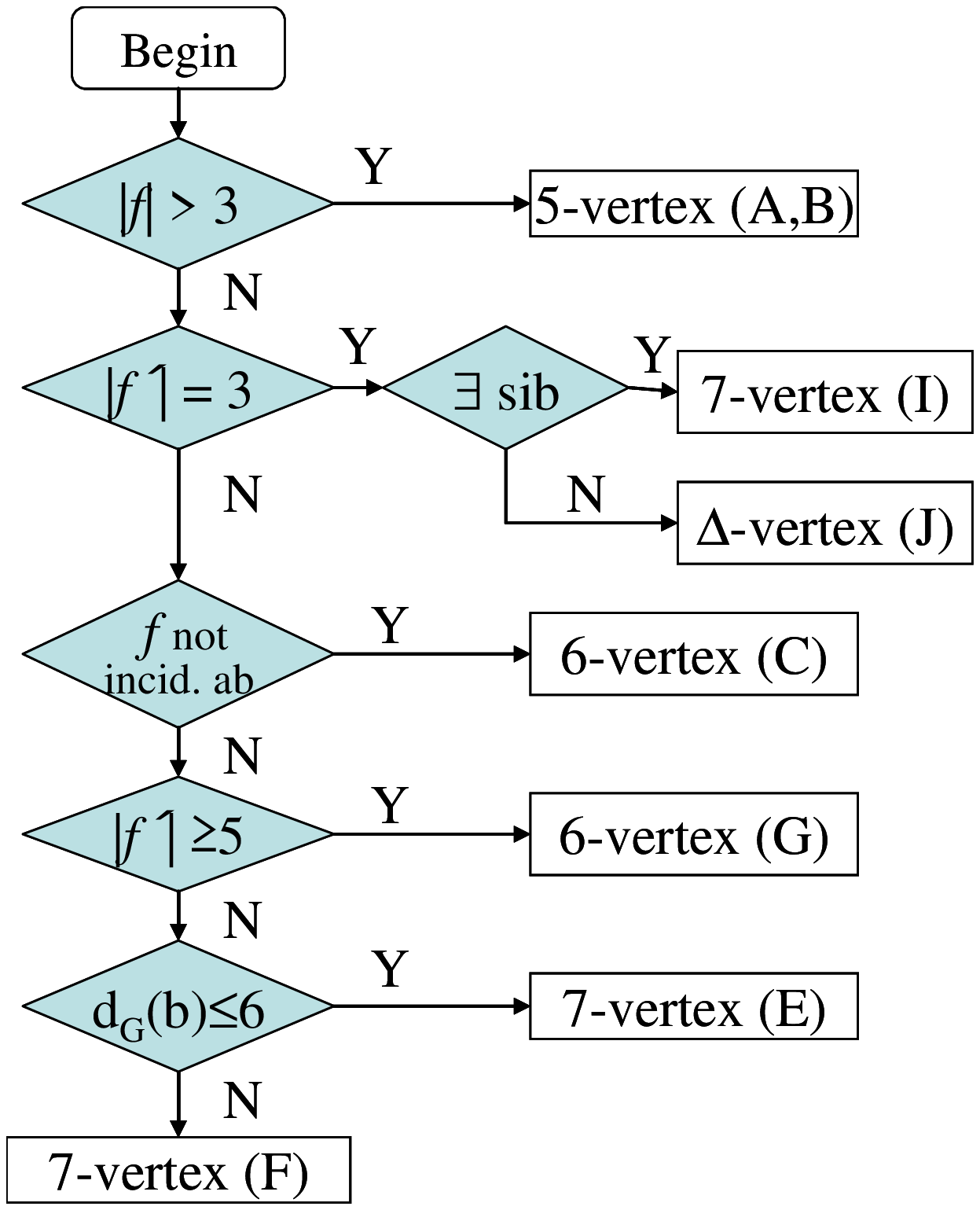}}
\caption{Reduction flowcharts}
\label{fig:flowcharts} 
\end{figure}

\begin{figure}[htbp]
\subfigure[$|f| \ge 5 \Rightarrow$ 4-vertex]{%
\label{subfig:reda}%
\begin{minipage}[b]{0.33\textwidth}
\centering 
\includegraphics[scale=\redscale]{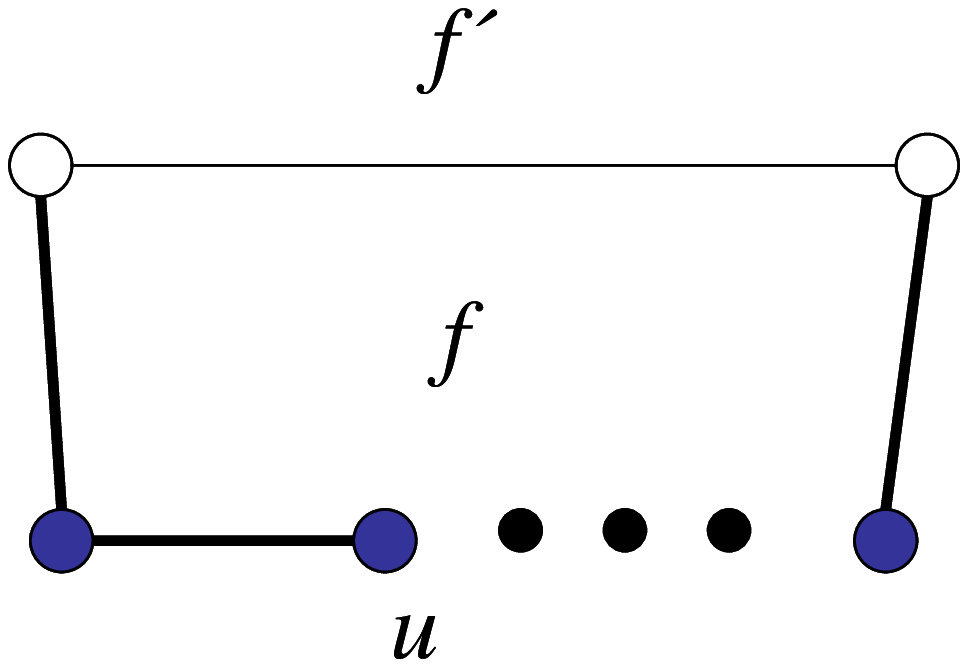}%
\end{minipage}%
}
\subfigure[$|f|=4 \Rightarrow$ 5-vertex]{%
\label{subfig:redd}%
\begin{minipage}[b]{0.33\textwidth}
\centering \includegraphics[scale=\redscale]{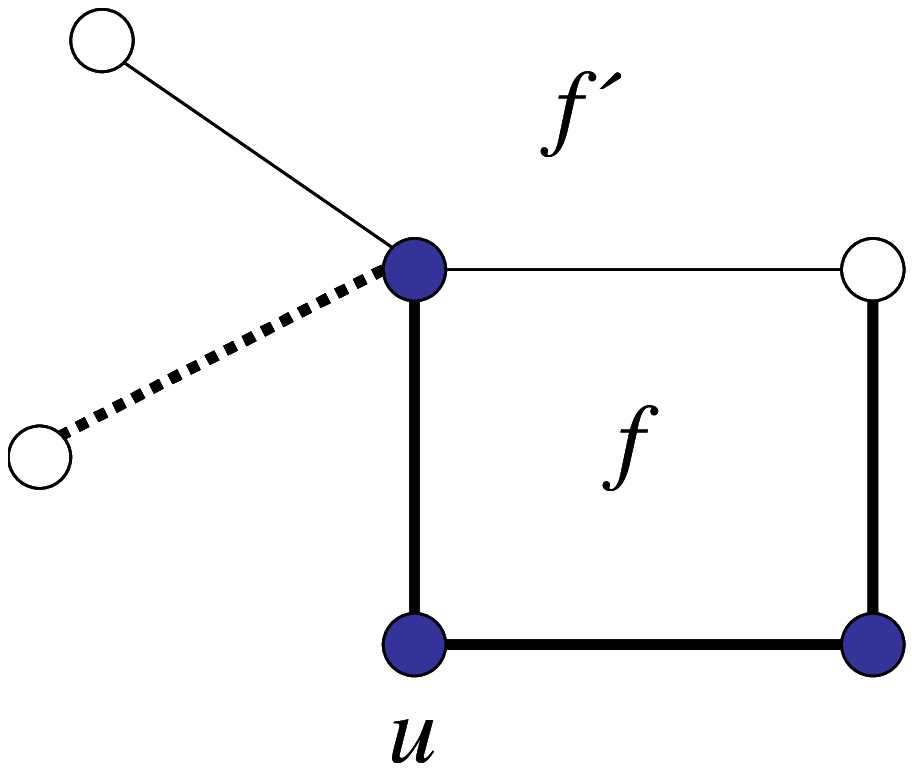}%
\end{minipage}%
}
\caption{Configurations for the case $|f|\ge 4$}
\label{fig:red1}
\end{figure}

\begin{figure}[htbp]
\addtocounter{subfigure}{2}
\subfigure[$|f|=3$, $|f'|\ge 4$, $f$ not incid.~on $ab$ $\Rightarrow$ 6-vertex]{%
\centering
\label{subfig:rede}%
\begin{minipage}[b]{0.33\textwidth}
\centering 
\includegraphics[scale=\redscale]{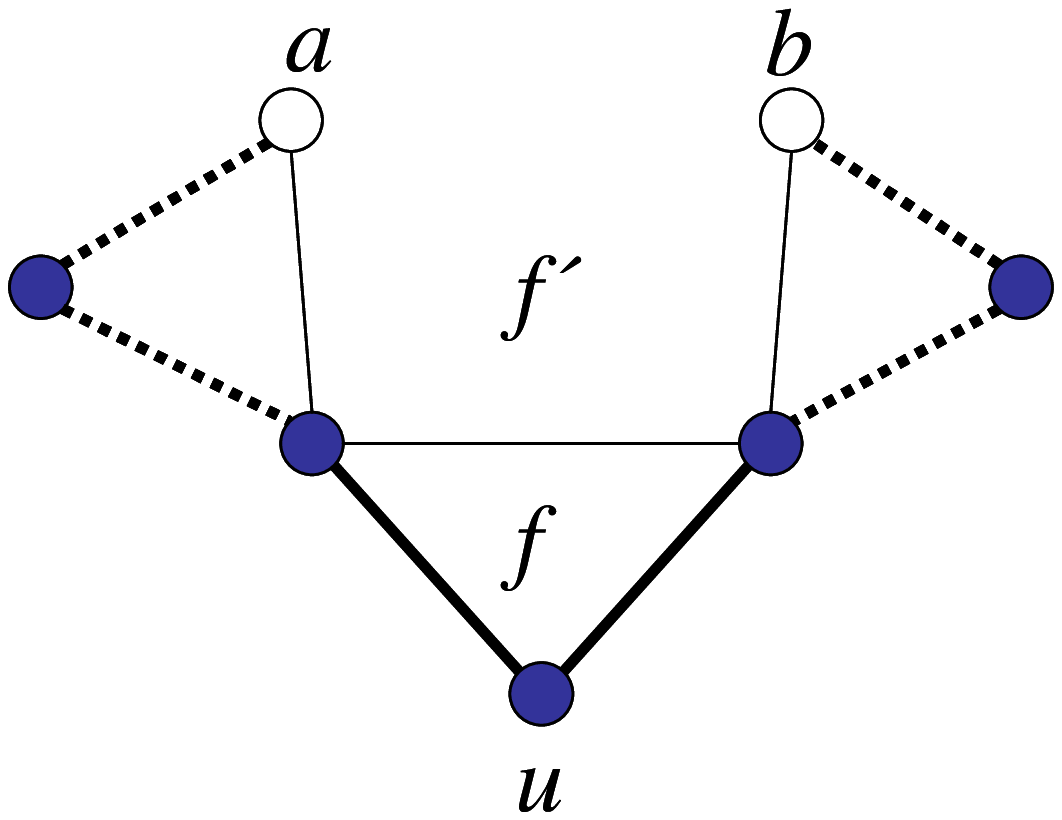}
\end{minipage}%
}
\subfigure[$|f|=3$, sib missing $\Rightarrow$ $(\Delta+1)$-vertex]{%
\label{subfig:redh}%
\begin{minipage}[b]{0.33\textwidth}
\centering 
\includegraphics[scale=\redscale]{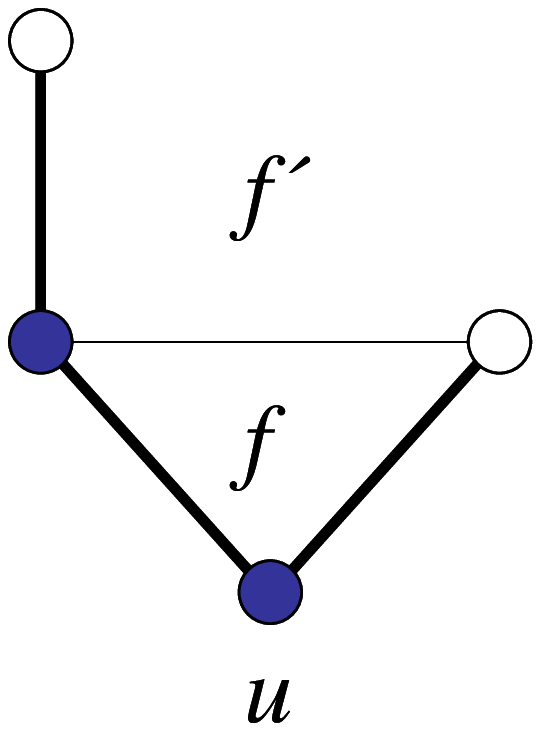}
\end{minipage}%
}
\subfigure[$|f|=3$, $|f'|=4$, all children of $f'$ incid.\ on $ab$, 
$d_G(b) \le 6$ $\Rightarrow$ $7$-vertex]{%
\label{subfig:redi}%
\begin{minipage}[b]{0.33\textwidth}
\centering 
\includegraphics[scale=\redscale]{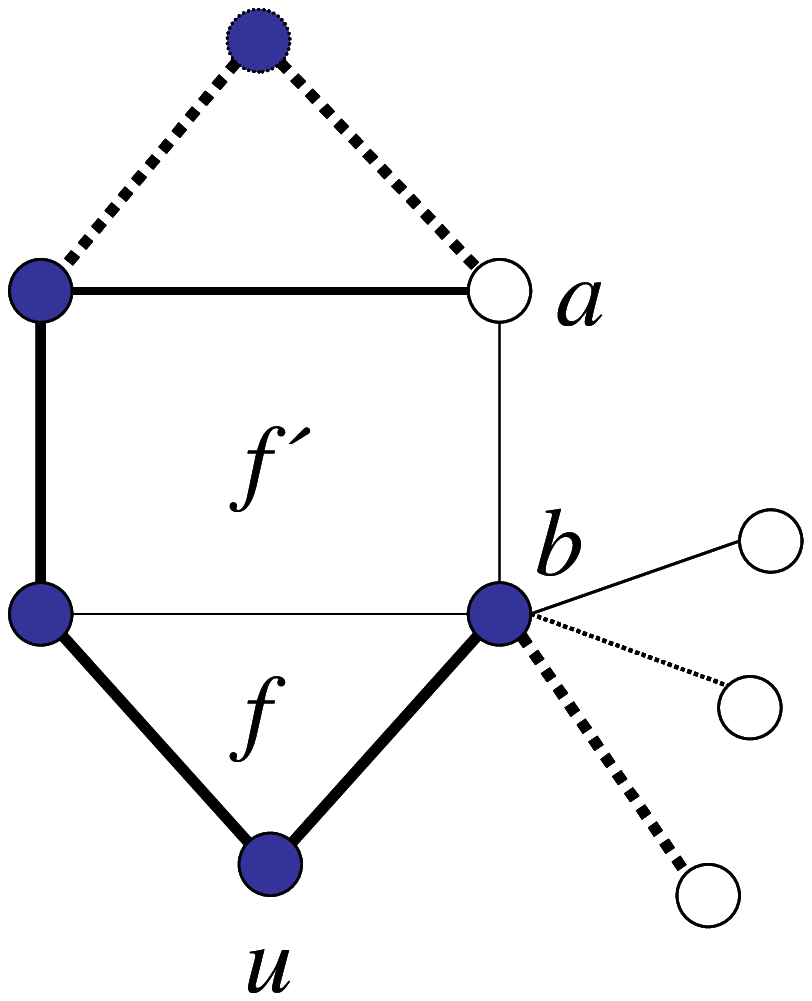}
\end{minipage}%
}
\subfigure[$|f|=3$, $|f'|=4$, no sib, $d_G(a)\le 5$  $\Rightarrow$ $7$-vertex]{%
\label{subfig:redj}%
\begin{minipage}[b]{0.33\textwidth}
\centering 
\includegraphics[scale=\redscale]{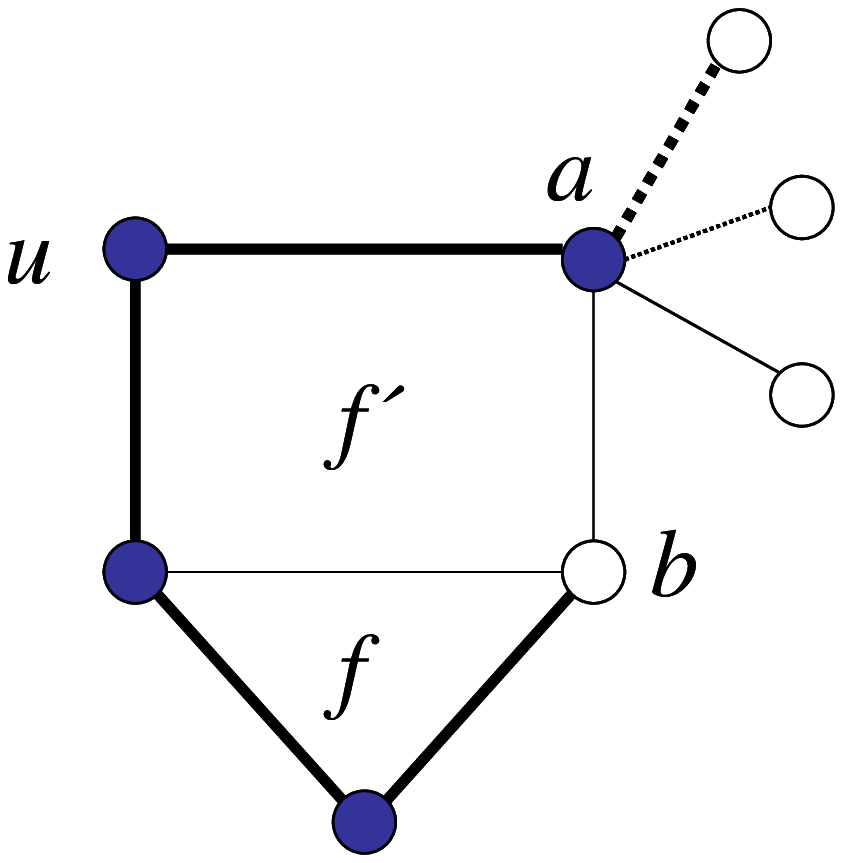}
\end{minipage}%
}
\subfigure[$|f|=3$, $|f'|\ge 5$, all children of $f'$ incid.\ on $ab$ 
  $\Rightarrow$ $6$-vertex]{%
\label{subfig:redg}%
\begin{minipage}[b]{0.33\textwidth}
\centering 
\includegraphics[scale=\redscale]{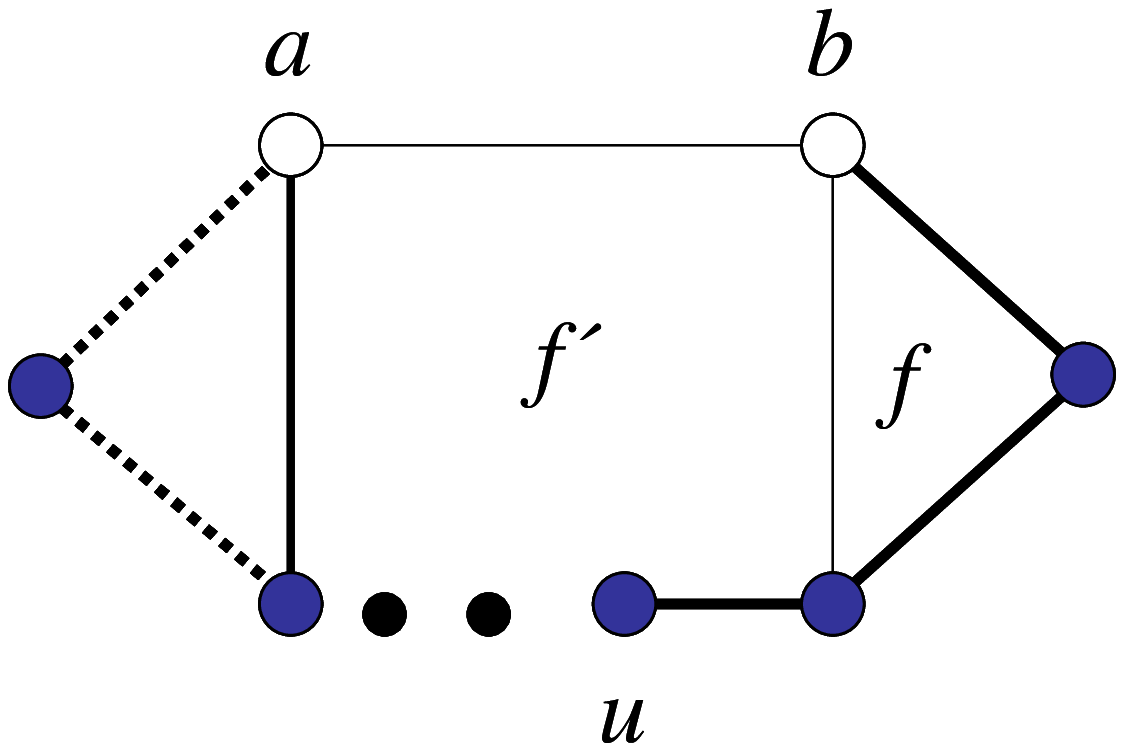}
\end{minipage}%
}
\caption{Configurations for the case $|f|=3$ and $|f'|\ge 4$}
\label{fig:red2}
\end{figure}

\begin{figure}[htbp]
\addtocounter{subfigure}{7}
\subfigure[$|f|=|f'|=3 \Rightarrow (\Delta+1)$-vertex]{%
\label{subfig:redk}
\begin{minipage}[b]{0.4\textwidth}
\centering \includegraphics[scale=\redscale]{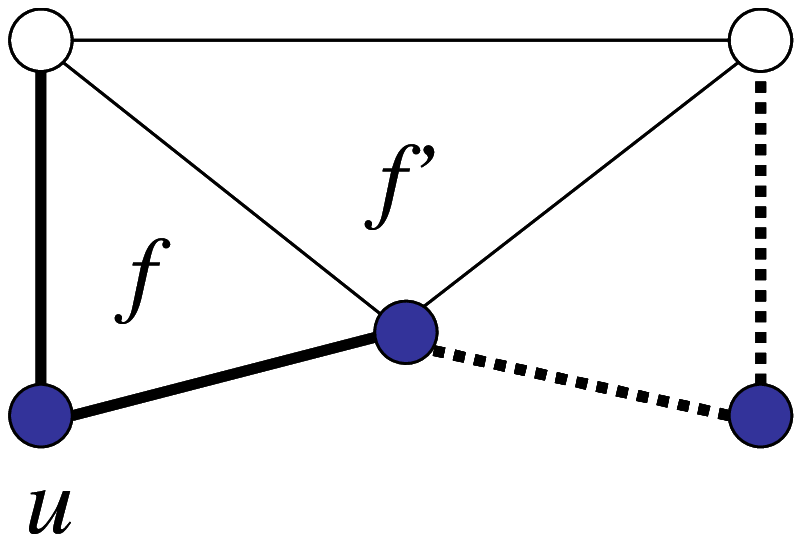}
\end{minipage}%
}
\subfigure[$|f|=|f'|=3$, $\nexists \mbox{sib} \Rightarrow \Delta$-vertex]{%
\label{subfig:redl}%
\begin{minipage}[b]{0.30\textwidth}
\centering \includegraphics[scale=\redscale]{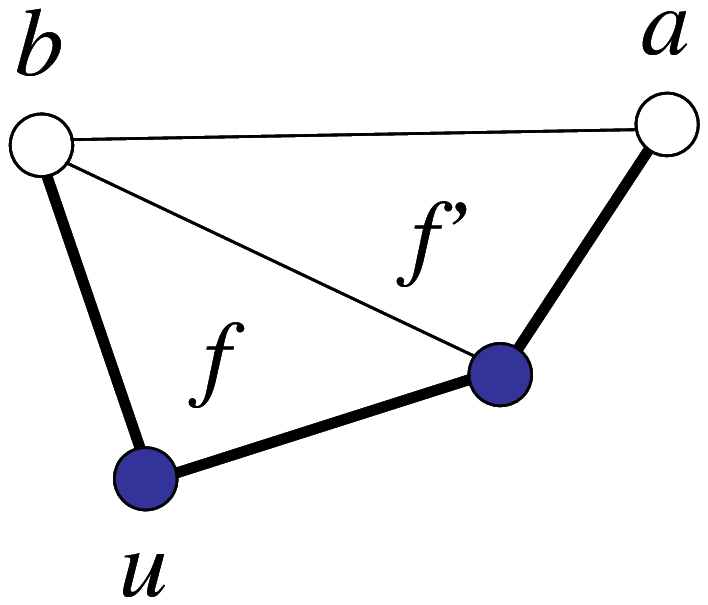}
\end{minipage}%
}
\subfigure[$|f|=|f'|=3, \mbox{$\exists$ sib} \Rightarrow 7$-vertex]{%
\label{subfig:redm}%
\begin{minipage}[b]{0.30\textwidth}
\centering \includegraphics[scale=\redscale]{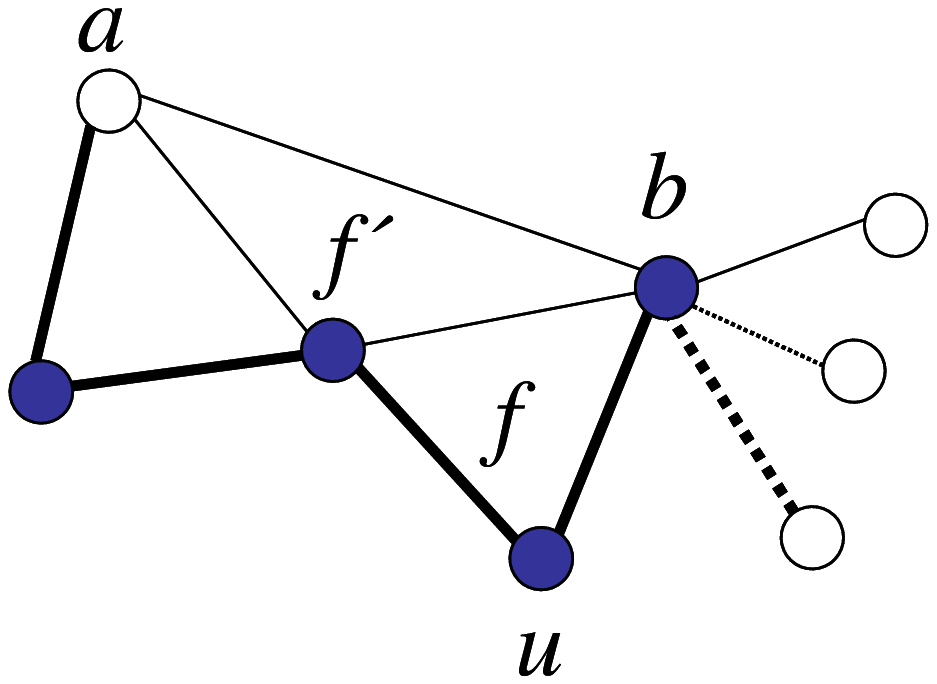}
\end{minipage}%
}
\caption{Configurations for the case $|f|=|f'|=3$}
\label{fig:red3} 
\end{figure}

We will give an exhaustive decision
tree, or a flowchart, that leads to a reducible configuration: 
A $(\Delta+1)$-reducible one when $\Delta\ge 5$ in 
Figure~\ref{fig:subfig:d1-flow},
and $\Delta$-reducible when $\Delta\ge 7$ in 
Figure~\ref{fig:subfig:d-flow}.
Each of the boxes of the branches 
corresponds to one of the reducible configurations 
of Figures~\ref{fig:red1}-\ref{fig:red3}, to be described shortly.

We shall assume that we are given $f$, $f'$, and $ab$ as promised by
Lemma \ref{lmm:good-face}. In the flowchart we use the following
notation: Recall that the cardinality of a face $f$, denoted $|f|$, is its number
of vertices.  $f$ is incident on $ab$ if it contains one of its
vertices.  We shall assume that if $f$ is incident on one of the
vertices of $ab$, then that vertex will be named $b$.
When traversing the flowchart (either one in Figure~\ref{fig:flowcharts}), we 
shall assume that all sibs of a face
are tested for a Y branch before proceeding to the corresponding N branch.

Each of the subfigures (A)-(J) in Figures~\ref{fig:red1}-\ref{fig:red3}
gives a configuration with a $k$-vertex marked as $u$.
Each of them expresses more generally a collection of configurations,
allowing for optional vertices as well as symmetric translations.
Edges that lie on the infinite face are shown in bold, while internal
edges are thin. Optional vertices and edges are shown with dotted edges.
Vertices shown in white have possible additional edges, while 
all neighbors of dark vertices (in blue) are shown in the figure.
We mark an $i$-ss face $f$ in the figure, along with its
parent $f'$.

\renewcommand{\thesubfigure}{\Alph{subfigure}}
\makeatletter
\renewcommand{\@thesubfigure}{(\thesubfigure)\space}
\renewcommand{\p@subfigure}{}
\makeatother
\subfigcapmargin = 10pt

\renewcommand{\thesubfigure}{\roman{subfigure}}

A set of configurations is \emph{unavoidable} for a class of graphs if 
every graph in the class contains at least one configuration from the set.
Our main technique, that bears a slight resemblance to the four color
theorem~\cite{AppelHaken}, is to give an unavoidable set of reducible configurations.

For reductions (B), (C), (D), (G), (H), and (I) to apply, $f$ must be 1-ss,
while for (E), (F), or (J) to apply, $f$ must be 2-ss.
\begin{lemma}
\label{lmm:Dleq5}
$\{A,B,C,D,H\}$ from Figures~\ref{fig:red1}-\ref{fig:red3}
is an unavoidable set of $(\Delta+1)$-configurations, 
for the class of outerplanar graphs of maximum degree $\Delta \ge 5$ 
with no simple leaf blocks.
\end{lemma}
\begin{proof}
We basically go through the flowchart in Figure~\ref{fig:subfig:d1-flow}.
Assume $f$, $f'$, and $ab$ as promised by Lemma \ref{lmm:good-face}.
If $|f|\ge 4$, then we have either case (\ref{subfig:reda}) or (\ref{subfig:redd}).
Assume then that $|f|=3$ and consider its parent face $f'$.
If $|f'|=3$, then we have case (\ref{subfig:redk}); otherwise, assume $|f'|\ge 4$.
If $f$ (or one of its sibs) is not incident on $ab$, 
then case (\ref{subfig:rede}) holds. Otherwise, $f$ is missing a sib on at 
least one side, so some edge of $f'$ borders the infinite face, in which case
(\ref{subfig:redh}) holds. In each case we obtain a $(\Delta+1)$-vertex 
and hence we have the lemma.
\end{proof}
From the above proof and Figure~\ref{fig:subfig:d1-flow} we obtain 
the following.
\begin{corollary}
\label{cor:Dleq5ss}
In each of the configurations from the unavoidable set $\{A,B,C,D,H\}$
in Lemma~\ref{lmm:Dleq5}, our $(\Delta+1)$-vertex 
is on the boundary of the 1-ss face $f$.
\end{corollary}
We have similarly the following for $\Delta \geq 7$.
\begin{lemma}
\label{lmm:Dleq7}
$\{A,B,C,E,F,G,I,J\}$ from Figures~\ref{fig:red1}-\ref{fig:red3}
is an unavoidable set of $\Delta$-configurations, 
for the class of outerplanar graphs of maximum degree $\Delta \ge 7$ 
with no simple leaf blocks.
\end{lemma}
\begin{proof}
We traverse the flowchart in Figure~\ref{fig:subfig:d-flow}.
Assume $f$, $f'$, and $ab$ as promised by Lemma \ref{lmm:good-face}.
If $|f|\ge 3$, then we have either of cases (\ref{subfig:reda}) 
and (\ref{subfig:redd}). Hence, we assume from now that $|f|=3$.

Consider the case $|f'|=3$. If $f$ has no sibling, 
then the case (\ref{subfig:redl}) applies and we have a $\Delta$-vertex.
Otherwise, $f$ has a sibling, which we can (by the above) assume
is also a 3-face. Since both $f$ and its sibling are 2-ss, then
one of them is bounded by three vertices of degree 2, 4 and at most 6 in $G$.
W.l.o.g.~we may assume $f$ to be this very face, in which case  
(\ref{subfig:redm}) applies and we have a $7$-vertex.

Consider now the case $|f'|\geq 4$. If $f$ is not incident on $ab$,
then the case (\ref{subfig:rede}) applies and we have a $6$-vertex.
Otherwise, assume $f$ (and all of its possible sibs) is incident on $ab$, 
in particular on $b$.
If $|f'| \ge 5$, then $f'$ has a 6-vertex 
as indicated in case (\ref{subfig:redg}); otherwise, assume $|f'|=4$.
Note that we are under the assumption that $f$ has no adjoint sibling
(since in that case we would have (\ref{subfig:rede})). 

By Claim~\ref{clm:ab-cd} we have that since $f$ is a 2-ss face, then 
$a$ and $b$ cannot both be of degree $\Delta\ge 7$.
Namely, only one of them can be incident on faces that descend from
$f'$ or its parent $f''$ (if it exists). 
The other has 3 edges incident on $f'$ and $f$ together, and at
most 3 edges incident on a sib of $f'$ and its possible child.

If $b$ (which is incident on $f$) has degree $6$ or less, then 
we have case (\ref{subfig:redi}), so $f$ contains a 7-vertex. Otherwise,
$a$ has degree $6$ or less.
We may for symmetric reasons assume that $f$ is the only child of $f'$
(as otherwise, the other child would work in the 
case (\ref{subfig:redi}) instead of $f$). 
Then, in fact, $a$ must have degree 5 or less, because it has only 2
edges incident on $f'$ and its children.
Then, the parent $f'$ contains a 7-vertex as indicated
in case (\ref{subfig:redj}), since the neighbors of the
unique degree-2 vertex on $f'$ have degree 3 and 5. 
This shows that in each case there is a
$\Delta$-vertex in $G$ and we have the lemma.
\end{proof}
Unlike the previous case of $\Delta = 5$, it is not always the case 
that the face $f$ contains
a $\Delta$-vertex when $\Delta\geq 7$. By 
Lemmas~\ref{lmm:Dleq5} and~\ref{lmm:Dleq7} we have
proved Theorem~\ref{thm:main-result} in the cases when 
$\Delta(G) \ge 5$.

We complete the proof by finishing the low-degree cases.


\paragraph{Cases with $\Delta \le 4$.}
For $\Delta=2$ we have $\ind(G^2)\leq \Delta + 2$.
In fact we have 
\begin{eqnarray*}
 \ind(G^2) & = & \left\{ \begin{array}{ll} 2 & \mbox{ for }G = P_k, \ \ k\geq 3,
                                               \mbox{ and }G = C_3, \\
                                           3 & \mbox{ for }G = C_4, \\
                                           4 & \mbox{ for }G = C_k, \ \ k \ge 5. 
                         \end{array} 
                 \right.
\end{eqnarray*}

For $\Delta \in \{3,4\}$ we have the following.
\begin{lemma}
\label{lem:ind-34}
For a outerplanar graph $G$ with $\Delta = k\in\{3,4\}$,
we have $\ind(G^2) \le 2k-2$.
\end{lemma}
\begin{proof}
By Lemma~\ref{lmm:ind-leq-k}, it suffices to show 
that $G$ contains a $(2\Delta-2)$-vertex.
If $G$ contains a degree-1 vertex, then it is a $\Delta$-vertex.
Otherwise, let $f$ be a leaf face in the dual tree $T^*(G)$.
If $|f|\ge 5$, then $f$ has a 4-vertex,
while if $|f|= 4$ then either of the degree-2 vertices of $f$ are $\Delta+1$-vertices.
Finally, if $|f|=3$, then the two neighbors of the degree-2 vertex
$u$ have at most $\Delta-2$ additional neighbors each. 
Hence, $u$ is a $2\Delta-2$ vertex. 
\end{proof}


\subsection{Choosability and algorithmic concerns.}
\label{sec:choosability}
As mentioned in Section~\ref{sec:prelim}, the bound on the inductiveness
of Theorem~\ref{thm:main-result} implies that {\Greedy} finds 
an optimal coloring of squares of outerplanar graphs of degree $\Delta \ge 7$.
When $\Delta \leq 6$, we can also obtain an efficient time
algorithm from the observation of Krumke, Marathe and Ravi~\cite{KMR:dialm}
that squares of outerplanar graphs have treewidth at most $3\Delta-1$.
This allows for the use of $2^{O(k)} n$-time algorithm for coloring
graphs of treewidth $k$. 
\begin{theorem}
\label{thm:linear-coloring}
There is a linear time algorithm to color squares of outerplanar graphs.
\end{theorem}

\paragraph{List coloring.}
Our approach for coloring $G^2$ for an outerplanar graph $G$ also
yields results regarding 
the list coloring, a.~k.~a.~{\em choosability}, of $G^2$ as well.
\begin{definition}
\label{def:choose}
A graph $G$ is \defn{{\em $k$-choosable}}
if for every collection of lists
$\{ S_v : v\in V(G)\}$ of colors where $|S_v| = k$ for 
each $v\in V(G)$, there is a coloring
$c : V(G) \rightarrow \bigcup_{v\in V(G)}S_v$,
such that $c(v)\in S_v$ for each $v\in V(G)$.
The minimum such $k$ is called the {\em choosability} or
the {\em list-chromatic number} of $G$,
and is denoted by \mathdefn{\ch(G)}.
\end{definition}
Note that if a graph is $k$-choosable, then it is $k$-colorable.
Also, by an easy induction, 
we see that if a graph is $k$-inductive then it is $(k+1)$-choosable.
For any graph $G$ we therefore have  $\chi(G)\leq \ch(G)\leq \ind(G) + 1$.

We thus obtain the following bound on choosability.
\begin{corollary}
\label{cor:gen-choose}
For any outerplanar graph $G$ with maximum degree $\Delta\geq 7$,
we have $\ch(G^2)=\Delta + 1$ and this is optimal.
\end{corollary}


\section{Chordal outerplanar graphs}
\label{sec:chordal}

Before we consider in detail the clique number and the chromatic
number for $G^2$ for an outerplanar graph in general, we will deal
with the chordal case first. This is because many chordal examples will
provide the matching lower bounds for the inductiveness, clique number
and the chromatic number as well. Here in the chordal case
we are able to present some structural results of $G$ in addition
to tight bounds of the three coloring parameters.

\begin{figure}[htbp]
\centering
\subfigure[$F_4 = \widehat{K_3}$]{
\includegraphics[scale=0.4]{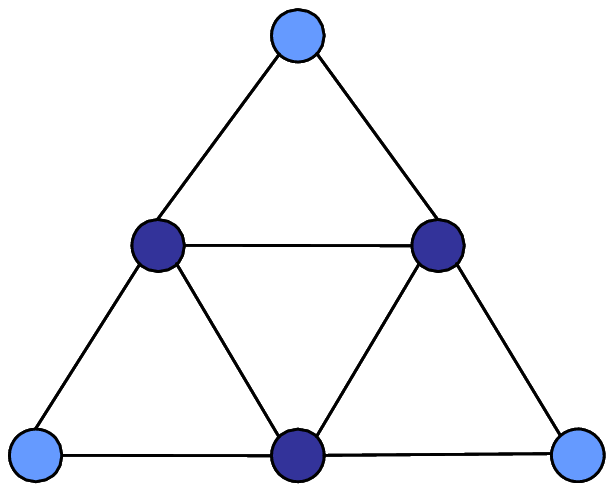}
\label{fig:subfig:f4}
}
\hspace{1cm}
\subfigure[$F_5 = \widehat{RL_4}$]{
\includegraphics[scale=0.4]{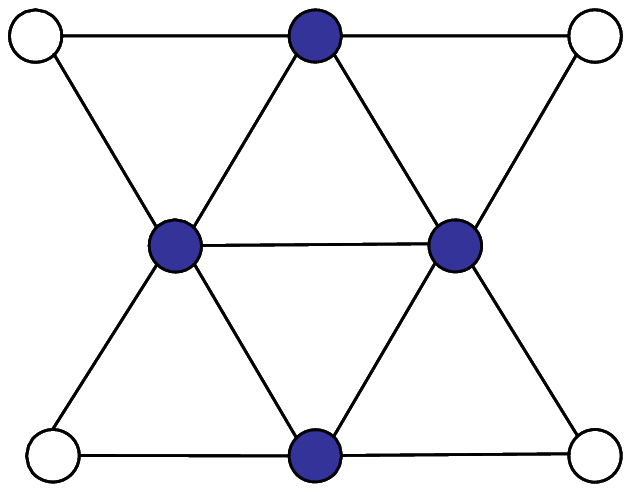}
\label{fig:subfig:f5}
}
\hspace{1cm}
\subfigure[$F_6 = \widehat{\widehat{K_3}}$]{
\includegraphics[scale=0.4]{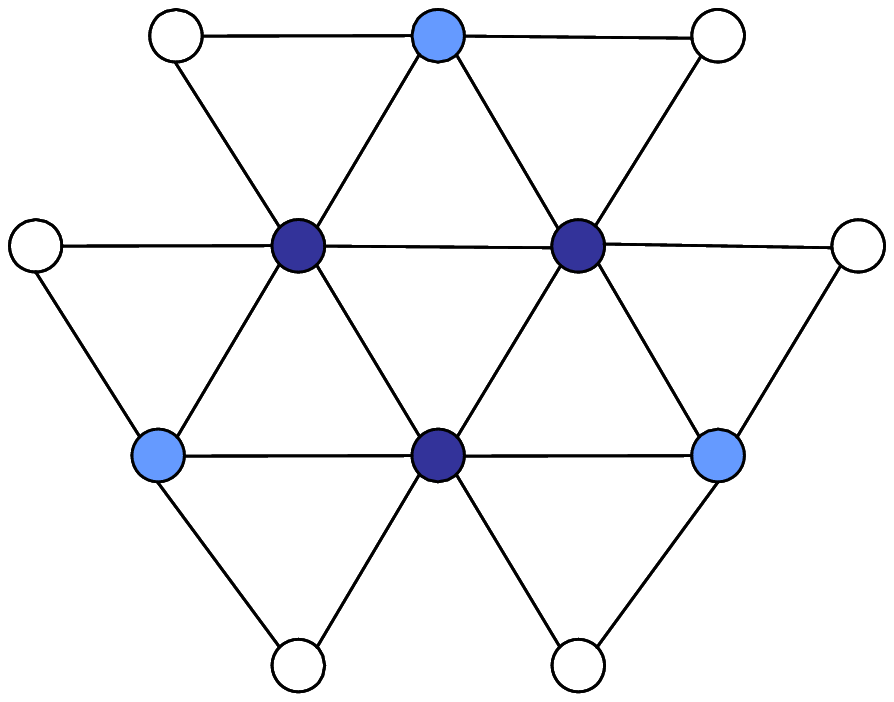}
\label{fig:subfig:f6}
}
\caption{Chordal outerplanar graphs of degree $\Delta=4, 5, 6$, with $\ind = \Delta+1$.}
\label{fig:f456}
\end{figure}

{\bf Conventions:} (i) Let $G$ be a given biconnected outerplanar 
on $n$ vertices of maximum degree $\Delta$,
with a fixed planar embedding. The graph obtain from $G$ by 
connecting an additional vertex to each pair of endvertices
of an edge bounding the infinite face, will be denoted
by $\widehat{G}$. Clearly $\widehat{G}$ will be an outerplanar graph on
$2n$ vertices of maximum degree $\Delta +2$.
(ii) By the {\em rigid $n$-ladder} or just the {\em rigid ladder} 
$RL_n$ on $n=2k$ vertices we will mean the graph given by 
\begin{eqnarray*}
V(RL_n) & =    & \{u_1,\ldots, u_k\}\cup \{v_1,\ldots, v_k\},\\
E(RL_n) & =    & \{\{u_i,v_i\}, \{u_i,u_{i+1}\}, \{u_i,v_{i+1}\}, \{v_i,v_{i+1}\} : \\
         &      & i\in \{1,\ldots,n-1\}\}\cup \{\{u_k,v_k\}\}. 
\end{eqnarray*}
For odd $n$, the graph $RL_n$ will mean $RL_{n+1}-u_{(n+1)/2}$.
(iii) Let $F_4 = \widehat{K_3}$, $F_5 = \widehat{RL_4}$, 
and $F_6 = \widehat{\widehat{K_3}}$, see Figure~\ref{fig:f456}.

Recall that when discussing the clique number or the chromatic number, 
we can by Lemma~\ref{lmm:biconn} assume $G$ to be biconnected.
One of the main results of this section is the following theorem.
\begin{theorem}
\label{thm:chordal-col}
For a chordal outerplanar graph $G$,
\[
\omega(G^2) = \chi(G^2) = \left\{
\begin{array}{cl}
\Delta+2 & \mbox{ if $\Delta=4$ and $F_4\subseteq G$,} \\
\Delta+1 & \mbox{ in all other cases.}
\end{array}
\right.
\]
\end{theorem}
We also derive a similar characterization of their inductiveness.

First note that the case $\Delta = 2$ is trivial, since there is only
one biconnected chordal outerplanar graph, namely $G = K_3$.

The case $\Delta \le 3$ is easy, since there are only 
three biconnected chordal outerplanar graphs 
with $\Delta \leq 3$:  $RL_2 = P_2$ the 2-path,
$RL_3 = C_3=K_3$ the 3-cycle, and $RL_4$ is the 4-cycle with one diagonal. 
From this we deduce the following tree-like structure of $G$ in this case.
\begin{lemma}
\label{lmm:D=3}
Let $G$ be a chordal outerplanar graph of maximum degree
$\Delta \le 3$. Then the blocks of $G$ are among
$\{RL_2,RL_3,RL_4\}$, where any two blocks from 
$\{RL_3,RL_4\}$ are separated by at least one $RL_2$ block.
\end{lemma}
Considering the leaf blocks of $G$, we obtain from the structure  
given in Lemma~\ref{lmm:D=3} the following.
\begin{theorem}
\label{thm:D=2,3}
For a chordal outerplanar graph $G$ with $\Delta\in\{2,3\}$,
we have 
\[
\omega(G^2) = \chi(G^2) = \ind(G^2)+1 = \Delta + 1.
\]
\end{theorem}

The case $\Delta = 4$ is more interesting, since it is the first
case involving a ``forbidden subgraph'' condition for both 
the clique and the chromatic number of $G^2$.
By considering the removal of a degree-2 vertex from $G$, we obtain
the following by induction on $n = |V(G)|$.
\begin{lemma}
\label{lmm:D=4}
A graph $G$ is a biconnected chordal outerplanar graph with $\Delta = 4$ if,
and only if, $G\in \{F_4\}\cup \{RL_n : n\geq 5\}$.
\end{lemma}
\begin{proof}
Clearly each graph in $\{F_4\}\cup \{RL_n : n\geq 5\}$ is biconnected
and outerplanar. Conversely, let $G$ be a biconnected outerplanar graph
on $n\geq 5$ vertices with maximum degree four. 
By removing a vertex $u$ of degree 2 from $G$, we obtain 
a biconnected outerplanar graph $G-u$ with $\Delta(G-u)\in \{3,4\}$,
and hence equal to $RL_4$ or, by induction, from 
the set $\{F_4\}\cup \{RL_n : n\geq 5\}$. 
Since $G$ is of maximum degree $\Delta = 4$ it is impossible that
$G-u = F_4$. For the same reason if $G-u = RL_4$, then $G = RL_5$. 
Also, $G-u = RL_5$ only when $G\in \{RL_6, F_4\}$,  and lastly if $G-u = RL_n$
for some $n\geq 6$, then $G = RL_{n+1}$ must hold, 
thereby proving the lemma.
\end{proof}
Note that for $G = F_4$ we have $G^2 = K_6$. Hence, in this case 
$\omega(G^2) = \chi(G^2) = 6 = \Delta + 2$, while $\ind(G^2) = 5$.

Observe that $\ind(RL_n^2) = 4$, for any $n\ge 5$, since removing the
last vertex in the square graph leaves the graph $RL^2_{n-1}$.
Thus, $\omega(RL_n^2)=\chi(RL_n^2) = 5$.
By Lemmas~\ref{lmm:biconn} and~\ref{lmm:D=4} we have the following.
\begin{theorem}
\label{thm:D=4}
For a chordal outerplanar graph $G$ with $\Delta=4$, we
have 
\[
\omega(G^2) = \chi(G^2) = \ind(G^2)+1 = 
\left\{
\begin{array}{cl}
6 & \mbox{ if }F_4\subseteq G, \\
5 & \mbox{ otherwise.}
\end{array}
\right.
\]
\end{theorem}

So far we have characterized chordal outerplanar graphs in terms
of the clique number, chromatic number and inductiveness of their squares
when $\Delta\in\{2,3,4\}$. 
Before we continue with the analysis of the chordal cases 
of $\Delta\in\{5,6\}$, we need the following definition and a lemma.
\begin{definition}
\label{def:separator}
Let $G$ be a graph. 
Call a subgraph $H\subseteq G$ on $h$ vertices 
an {\em $h$-separator}, or just a {\em separator} if
it induces a clique in $G^2$ whose removal breaks $G^2$ into
disconnected components.
\end{definition}
The following lemma shows that it suffices to bound the clique number and
chromatic number for graphs without separators.
\begin{lemma}
\label{lmm:chrom-sep}
Let $G$ be a graph and $H$ a separator of $G$  
with $G = G'\cup G''$ and $H = G'\cap G''$. Then we have 
\begin{eqnarray*}
\omega(G^2) & = & \max\{\omega({G'}^2),\omega({G''}^2)\}, \\
\chi(G^2)   & = & \max\{\chi({G'}^2), \chi({G''}^2)\}.
\end{eqnarray*}
\end{lemma}
\begin{proof}
Since $G^2 = {G'}^2\cup {G''}^2$ and ${G'}^2 \cap {G''}^2 = H^2$, which
is a clique, we have the stated clique number for $G^2$.
Further, any optimal coloring of either ${G'}^2$ or ${G''}^2$ 
can be extended to an optimal coloring of $G^2$ by a suitable permutation
of the colors. Hence we have the lemma.
\end{proof}
Recall that a tree is full if there are no vertices of 
degree 2.
\begin{lemma}
\label{lmm:full-D56}
Let $G$ be a biconnected chordal outerplanar graph with a
full weak dual $T^*(G)$.
\begin{enumerate}
  \item If $\Delta = 5$, then $G = F_5$.
  \item If $\Delta = 6$, then $G \in \{F_6\}\cup \{\widehat{RL_n} : n\geq 5\}$.
\end{enumerate}
\end{lemma}
\begin{proof}
Since $G$ is chordal it is uniquely determined by a plane 
embedding of its full weak dual $T^*(G)$. 
If $\Delta \in \{5,6\}$, then $T^*(G)$ contains
a path of length $\Delta - 2$ and since $T^*(G)$ is full
each internal vertex of this path must have degree of three.
Therefore $T^*(G)$ has at least $2(\Delta - 2)$ vertices
and so $G$ has $n\geq 2(\Delta - 1)$ vertices. 
The unique full tree on $2(\Delta - 2)$ vertices is the weak dual of 
$G = \widehat{RL_{\Delta - 1}}$ on $2(\Delta - 1)$ 
vertices and with maximum degree $\Delta$.
Proceeding by induction, assume $G$ has $n\geq 2(\Delta - 1)$ vertices
has a full weak dual $T^*(G)$ on $n - 2\geq 2(\Delta - 2) $ vertices. 
Removing two siblings from $T^*(G)$
whose parent is a leaf in the pruned tree $\pr(T^*(G))$ corresponds
to removing two degree-2 vertices $u$ and $v$ from $G$ of distance 2 from each
other in $G$, and obtaining $G'' = G - \{u,v\}$. 
Now $\pr(T^*(G)) = T^*(G'')$ is full and has $n - 4\geq 2(\Delta - 3)\geq 4$ 
vertices, since its maximum degree is three.

If $T^*(G'')$ has four vertices, then it is the unique 4-star
and $G'' = F_4$. Hence $G = F_5$ and $\Delta = 5$ here.
Note that the number of vertices in any full tree with maximum
degree of three is always even. Hence, there are no full trees
on five or seven vertices.

If $T^*(G'')$ has six vertices, then $T^*(G'')$ is unique and
$G'' = F_5$. Hence $G = \widehat{RL_5}$ in this case and $\Delta = 6$.

Otherwise $T^*(G'')$ must have at least eight vertices and
hence $G''$ has at least ten vertices. By induction hypothesis we
have $G'' \in \{F_6\}\cup \{\widehat{RL_n} : n\geq 5\}$. To have
$G'' = F_6$ is impossible, since that would create $G$ with $\Delta = 8$.
So $G'' = \widehat{RL_n}$ for some $n\geq 5$.
If $G'' = \widehat{RL_5}$ then either $G = F_6$ or $G = \widehat{RL_6}$.
Otherwise $G = \widehat{RL_{n-1}}$, which completes the proof.
\end{proof}
\begin{theorem}
\label{thm:separator-D56}
Let $G$ be a biconnected chordal outerplanar graph with $\Delta \in\{5,6\}$.
If $T^*(G)$ is not full, then either $G^2$ is a clique or $G$  
has an $h$-separator where $h\in\{4,5,6,7\}$.
\end{theorem}
\begin{proof} 
Note that since $G$ is chordal, then every face corresponding to a vertex in
$T^*(G)$ is bounded by a triangle. So for a given $T^*(G)$,
the structure of $G$ is determined except for the degree-2 vertices of $G$.
By assumption $T^*(G)$ has a degree-2 vertex. We consider the
following three cases.

If $T^*(G)$ has no vertex of degree 3, then $T^*(G)$ is a simple path.
In this case it is easy to see that $G^2$ has a path decomposition consisting
of cliques, where each clique is induced by a closed neighborhood of a vertex
of $G$, necessarily of of size 5, 6 or 7. In particular, $G^2$ is 
by~\cite[Prop.~12.3.8]{Diestel} a chordal graph and is therefore a clique
or has a 4-separator.

Consider next the case where $T^*(G)$ has a degree-2 vertex $u^*$ that lies on 
a path connecting two degree-3 vertices of $T^*(G)$. If $e$ is the unique edge
in $G$ bounding the triangular face $f$ corresponding to $u^*$ and the infinite face, then
let $w\in V(G)$ be the vertex opposite the edge $e$ in the triangle $f$. Since $e$ and two
other edges $e'$ and $e''$ incident to $w$ in $G$ all bound the infinite face of $G$, 
we see that the closed neighborhood of $N_G[w]$ is an $h$-separator of $G$,
where $h = d_G(w) + 1\in \{5,6,7\}$.

Lastly, consider the case where every degree-2 vertex of $T^*(G)$ lies
on a path connecting a degree-3 vertex and a leaf of $T^*(G)$.
In this case $T^*(G)$ must contain a degree-2 vertex $v^*$ that has 
a leaf $u^*$ as a neighbor. Let $u,v,w\in V(G)$ be the vertices 
bounding the triangular face $f$ corresponding to the leaf $u^*$, where $u$
has degree 2 in $G$ and $v$ has degree 2 in $G - u$. Note that
$v^*$ is a leaf in $\pr(T^*(G))$. In this case $V = N_G[w]\setminus \{u\}$
induces a clique in $G^2$ which separates $u$ from the rest of the graph
$G^2$. Since $T^*(G)$ has a vertex of degree 3, then $G^2 - V$ has at least
two components where one component consists of the singleton $u$.
Hence, $V$ induces an $h$-separator of size $h = |V| = d_G(w)\in \{4,5,6\}$.
This completes our proof.
\end{proof}
Note that $F_5^2$ is $K_8$ with
two perpendicular diagonals removed when the vertices are located
on a regular 8-gon. Therefore these two pairs of opposite nonadjacent
vertices can be colored by the same color in $F_5^2$, and hence
$\chi(F_5^2) = \Delta + 1 = 6$. 

We note further that $F_6$ is biconnected chordal outerplanar graph with 
$\Delta = 6$. Also, the subgraph $\widehat{K_3}$ in $F_6$ induces a 
clique in $F_6^2$, and hence
each vertex there must have a unique color, say cyclically with colors
$1,2,3,4,5,6$ starting with a degree-2 vertex of $\widehat{K_3}$.
Of the remaining six vertices of $F_6$, color three of them with a new
color 7, all of distance three apart, and the remaining three with the
colors 1, 3 and 5. Hence $\chi(F_6^2) = \Delta + 1 = 7$.

To color ${\widehat{RL_n}}^2$ we can start by coloring the degree-2 vertices of
$\widehat{RL_n}$ alternatively with colors 1 and 2 cyclically. The rest of
the vertices, that constitute $RL_n$, we can then color with the remaining five
available colors, since we have already that $\ind(RL_n^2) = 4$.
Hence we have $\chi(RL_n^2) = \Delta(RL_n) + 1 = 7$. 

With the above in mind together with 
Lemmas~\ref{lmm:biconn}, \ref{lmm:chrom-sep}, and \ref{lmm:full-D56} 
and Theorems~\ref{thm:main-result} and~\ref{thm:separator-D56} 
we obtain in particular the following.
\begin{corollary}
\label{cor:chordal-D56}
For a chordal outerplanar graph $G$ with $\Delta\in\{5,6\}$ we
have 
\[
\omega(G^2) = \chi(G^2) = \Delta + 1.
\]
\end{corollary}
Corollary~\ref{cor:chordal-D56} together with Theorems~\ref{thm:D=2,3}
and~\ref{thm:D=4} complete the proof of Theorem~\ref{thm:chordal-col} as
well as the entries in Table~\ref{tab:Delta} in the 
chordal case for $\Delta\in\{2,3,4,5,6\}$.
\begin{observation}
\label{obs:D=456}
For each $\Delta\in \{4,5,6\}$, there are infinitely
many biconnected chordal outerplanar
graphs $G$ of maximum degree $\Delta$
with $\ind(G^2) = \Delta + 1$.
\end{observation}
\begin{proof}
It suffices to show that for each $\Delta \in \{4, 5,6\}$ there are 
infinitely many biconnected outerplanar graphs $G$ 
whose squares are of minimum degree $\Delta + 1$.
Refer to Figure~\ref{fig:f456} for the appearance of the graphs.

For $\Delta = 4$, consider $F_4 = \widehat{K_3}$,
whose square is $K_6$ and hence has a minimum degree $\Delta + 1$.
By fusing together edges whose endvertices have degrees
two and four, in two or more copies of $F_4$, we can
construct an infinite family of such biconnected outerplanar
graphs $G$ with $\Delta=4$ with $\ind(G^2) = 5$.

For $\Delta = 5$, consider $F_5$,
whose square is $K_8$ with two perpendicular diagonals
removed when the vertices are located on a regular 8-gon.
Also in this case we can fuse together two edges
with endvertices of degree 2 and 5, of two or more
copies of $F_5$, and obtain 
an infinite family of such biconnected outerplanar
graphs $G$ with $\Delta=5$ with $\ind(G^2) = 6$.

Finally, for $\Delta = 6$, consider $F_6$.
In this case we have $\ind(F_6^2) = 7$, and also 
here we can fuse edges with endvertices of degree
2 and 4, of two or more copies of $F_6$, 
to form an infinite family of biconnected 
outerplanar graphs $G$ of maximum degree $\Delta = 6$
and with $\ind(G^2) = 7$.
This completes the proof.
\end{proof}
{\sc Remark:} We note that $F_4^2 = K_6$ 
is chordal, but neither $F_5^2$ nor $F_6^2$ 
are chordal, showing that the square of a chordal
graph does not need to be chordal. 
This is consistent with the characterization of those chordal graphs
whose squares are chordal given
in~\cite{LaskarShier}.
\begin{theorem}
\label{thm:chordal-ind}
For a chordal outerplanar graph $G$, we have
$\ind(G^2) = \Delta$ or $\Delta+1$.
Necessary and sufficient conditions that $\ind(G^2) = \Delta+1$, are
one of the following: 
\begin{enumerate}
 \item $\Delta=4$ and $F_4 \subseteq G$,
 \item $\Delta=5$ and $F_5 \subseteq G$, or
 \item $\Delta=6$ and one of $\{ F_6\} \cup \{RL_n\}$, $n \ge 4$, is a
subgraph of $G$.
\end{enumerate}
\end{theorem} 


\section{Clique number}
\label{sec:clique}

In this section we deal with the clique number of $G^2$. This parameter
is the easiest to deal with. Nonetheless the exact computation of 
the clique number will rely on some results in previous
Sections~\ref{sec:ind} and~\ref{sec:chordal}. Recall that 
by Lemma~\ref{lmm:biconn} we can assume $G$ to be biconnected.
Further, $G$ is here not necessarily chordal.
\begin{lemma}
\label{lmm:D=5,6}
For an outerplanar graph $G$ we have $\omega(G^2) \le \Delta+2$, unless $G=C_5$.
If $\Delta\ge 5$, then $\omega(G^2) = \Delta+1$.
For $\Delta\ge 6$ we further have that any clique with 
$\Delta+1$ vertices is the closed neighborhood of some vertex.
\end{lemma}
\begin{proof}
Consider an induced subgraph $H_k$ of $G$ with $k+1$ vertices:
a vertex $u$ and its neighbors $u_1, u_2, \ldots, u_k$ in a clockwise order in the
plane embedding of $G$. Then only adjacent pairs $u_i$ and $u_{i+1}$, 
for $i=1,\ldots,k-1$, may be connected by an edge. Consider now a
vertex $w$ that is not a neighbor of $u$. Then, 
$w$ can be adjacent to at most two neighbors of $u$ and only
consecutive ones, by the outerplanarity property. 
If $k\ge 5$, then $w$ cannot be adjacent to both one of $u_1$ and $u_2$
and to one of $u_{k-1}$ and $u_k$. Thus, it must be of distance at
least 3 from either $u_1$ or $u_k$. Hence, $H_k \cup \{w\}$ is not a
clique in $G^2$. This shows that if a clique in $G^2$ contains
a closed neighborhood of a vertex of degree $k\geq 5$, then
the clique consists precisely of those $k+1$ vertices.

Consider now the case $k=4$ and we have two vertices $w_1$ and $w_2$
that are non-neighbors of $u$. In order to be of distance at most 2
from both $u_1$ and $u_4$, a vertex must be adjacent to $u_2$ and
$u_3$. But, in an outerplanar graph, not both $w_1$ and $w_2$ can be
so. Hence, $H_4 \cup \{w_1, w_2\}$ does not form a clique.
This shows that if a clique in $G^2$ contains a closed neighborhood
of a vertex of degree $k=4$, then the clique consists of at most
$k+2$ vertices, the vertices of the closed neighborhood plus
a possibly additional vertex. 

Suppose now an induced subgraph $H$ of maximum degree 3
induces a 6-clique in $G$.
Recall that by Lemma~\ref{lmm:biconn} we assume $H$ to be biconnected
and therefore induced by a cycle. There can be at most two chords in $H$ and 
they must be disjoint since
$\Delta(H) = 3$. Then there are two vertices in $H$ of degree 2 that 
are of distance 3 in $H$. Further, since all vertices of $H$ lie on the
outer face, there can be no vertex outside $H$ connecting them.

From the above paragraphs we conclude that if a clique of $G^2$ 
contains $\Delta+1\geq 7$ vertices, it must be a closed neighborhood
of a vertex of degree $\Delta$.
Hence, the lemma.
\end{proof}
We conclude this section by quickly discussing matching
lower bounds for $\omega(G^2)$.
Note, that we are here still under the assumption that $G$ is biconnected.

If $\Delta=2$, then $G = C_k$ is a cycle on $k\geq 3$ vertices
and we clearly have
\begin{eqnarray*}
 \omega(G^2) & = & \left\{ \begin{array}{ll} 
                              \Delta + 3 & \mbox{ if }G = C_5, \\
                              \Delta + 2 & \mbox{ if }G = C_4, \\
                              \Delta + 1 & \mbox{ otherwise.}  \\
                         \end{array} 
                 \right.
\end{eqnarray*}

For $\Delta=3$ the upper bound of $\Delta + 2$ is matched if $G$
is the graph obtained by adding one chord to the 5-cycle $C_5$.

For $\Delta\in\{4,5\}$ the matching upper bound of both $\Delta + 2$
when $\Delta = 4$ and $\Delta + 1$ when $\Delta=5$ is obtain
when $G = F_4$ shown in Figure~\ref{fig:f456}.

If $\Delta\geq 6$, then by the above Lemma~\ref{lmm:D=5,6},
the matching upper bound of $\Delta + 1$ is obtained
by a closed neighborhood of any vertex of degree $\Delta$.

Together with what was obtained in the previous 
Section~\ref{sec:chordal}, we therefore have all the 
entries for $\omega(G^2)$ for an outerplanar graph $G$ (chordal or not)
displayed in Table~\ref{tab:Delta}.


\section{The chromatic number}
\label{sec:col}

Recall that for an outerplanar graph $G$ with $\Delta\geq 7$
we have by Theorem~\ref{thm:main-result} that $\ind(G^2)\leq \Delta$
and hence $\chi(G^2)\leq \Delta +1$ which is optimal.
(In fact, Lemmas~\ref{lmm:biconn} and~\ref{lmm:Dleq7}
also imply this.) Hence, we will assume throughout this
section that $\Delta\leq 6$.


\subsection{Cases with $\Delta\in\{2,3,5\}$}

If $\Delta = 2$, then we have an upper bound 
$\chi(G^2) \leq \ind(G^2) + 1\leq \Delta + 3$
and a matching lower bound is obtained when
$G = C_5$. In fact, if $G$ is $P_k$ (resp.~$C_k$) the path (resp.~cycle) 
on $k$ vertices, then $\Delta = 2$ and 
it can be verified that for $k\geq 2$ we have that  
\begin{eqnarray*}
\chi(G^2) & = & \left\{ \begin{array}{ll} 
                          3 & \mbox{ if } G = P_k \cup C_{3k}, \\
                          4 & \mbox{ if } G = (C_{3k+1} \cup C_{3k+2}) \setminus C_5, \\
                          5 & \mbox{ if } G = C_5. 
                        \end{array} 
                \right. \\
\end{eqnarray*}
Further, {\Greedy} obtains an optimal coloring even when inductiveness 
is not a tight bound on the chromatic number, that is on $C_n$ for $n\ge 6$.

For $\Delta=3$ and $\Delta=5$, we have an upper bound of 
$\chi(G^2) \le \ind(G^2)+1 = \Delta + 2$.
If $G$ is the graph obtained by adding a chord to the
5-cycle $C_5$, then $G^2$ is a clique and hence
the lower bound for $\chi(G^2) = \Delta + 2$ is obtained
for $\Delta = 3$. 

\begin{figure}[htbp]
\centering
\subfigure[$G_{10}$, $\chi(G^2_{10})=7=\Delta+2$]{
\includegraphics[scale=0.4]{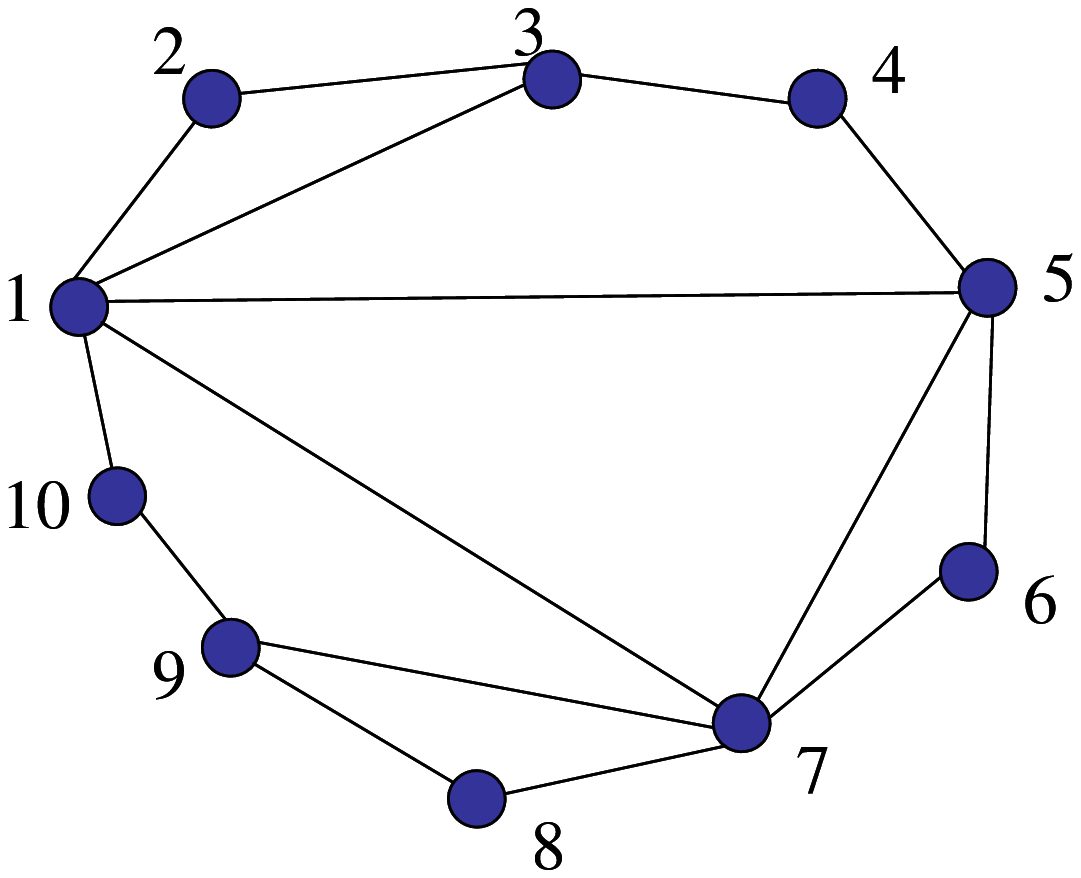}
\label{fig:subfig:g10}
}
\hspace{2cm}
\subfigure[$\overline{G^2_{10}}$, the complement of the square of $G_{10}$]{
\includegraphics[scale=0.4]{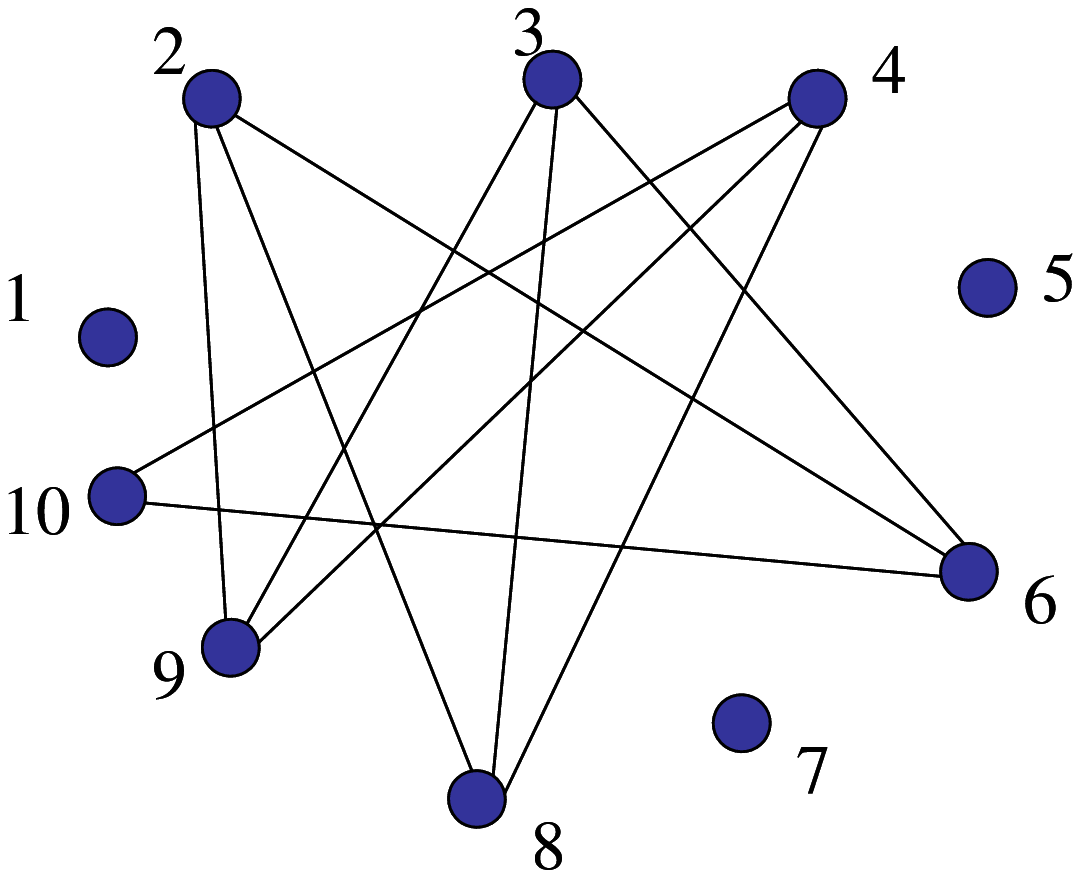}
\label{fig:subfig:g10-comp}
}
\label{fig:g10}
\caption{An example of a graph with $\Delta=5$ and $\chi = 7 = \Delta+2$}
\end{figure}

Consider now the case $\Delta = 5$. Let $G_{10}$ be the graph
on ten vertices given in (\ref{fig:subfig:g10}) in Figure~\ref{fig:g10}. 

To see that $G_{10}^2$ requires $7$ colors, it is easiest to try to
cover with cliques the complement graph $\overline{G_{10}^2}$, shown in
(\ref{fig:subfig:g10-comp}) in Figure~\ref{fig:g10}. 
Each of the vertices $u_1$, $u_5$ and $u_7$ require their own clique, 
while for the remaining 7 vertices, there is no 3-clique. Hence, 7 cliques are
required to cover $\overline{G_{10}^2}$. That is, 7 colors are required 
to color $G_{10}^2$. We note that $G_{10}$ has four edges with 
endvertices of degree 2 and 3 respectively.  By fusing together 
two copies of $G_{10}$ along
these edges in such a way that a degree-2 vertex in one copy is
identified with a degree-3 vertex in another copy, we can make an
infinite family of outerplanar graphs with $\Delta = 5$, such that their
square has chromatic number of $7$. We summarize in the following.
\begin{theorem}
\label{thm:D=5-infinite}
There are infinitely many biconnected outerplanar graphs $G$ with 
maximum degree $\Delta = 5$ such that $\chi(G^2) = 7$.
\end{theorem}


\subsection{Cases with $\Delta=4$}

Although it is impossible to obtain the tight 
bound of $\chi(G^2)$ in the case when $\Delta = 4$
via the inductiveness of $G^2$, our approach here
will in similar fashion be inductive: We show how we
can extend an optimal coloring of the square of
a subgraph $G'$ of $G$ to that of $G$.
\begin{lemma}
\label{lmm:chi-Delta=4}
If $G$ is an outerplanar graph with maximum degree $\Delta=4$, 
then $\chi(G^2)\leq 6$.
\end{lemma}
\begin{proof}
By Lemma~\ref{lmm:biconn} we may assume $G$ to be biconnected.
Since $\Delta = 4$, then $\pr(T^*(G))$ is a proper tree and
hence each 0-ss face $f$ has a parent $f'$. 
If we can find a $k$-vertex in $G$, for $k\le 5$, then a 6-coloring of
$G$ follows by induction on $|V(G)|$.
Observe that if $|f|\ge 4$,
then each degree-2 vertex of $f$ is a 5-vertex. If $|f|=3$ and $|f'|=3$,
then the unique degree-2 vertex on $f$ is a 5-vertex.
Additionally, if $|f|=3$ and $f'$ has a face bounding the infinite
face, then one 0-ss child of $f'$ has a 5-vertex.
Hence, we can assume that each 0-ss face is a 3-face
whose parent $f'$ satisfies $|f'|\geq 4$ and has no edge bounding the
infinite face. 

This implies one of the following two possibilities: (i) 
$f'$ has a grandparent $f''$, or (ii) every edge of $f'$ 
also bounds a child of $f'$. In this case we can simply
choose one designated child to play the role of the grandparent $f''$.

Let $v_1,\ldots,v_k$ be the vertices of $f'$ where $k\geq 4$, 
and $v_1v_k$ be the edge bordering $f''$. Let $u_1,\ldots,u_{k-1}$ be the vertices such that
$v_iu_iv_{i+1}$ is a 0-ss child of $f'$, for $i=1,\ldots, k-1$.
Finally, since $\Delta=4$, $v_1$ and $v_k$ have at most one additional
neighbor each; denote them by $x$ and $y$, respectively, which are not
necessarily distinct.

Let $G'$ be the graph obtained after eliminating the face $f'$ and its
children (i.e., removing vertices 
$v_2, \ldots, v_{k-1}$, $u_1, \ldots, u_{k-1}$ and incident edges from $G$.)
We show that a 6-coloring of ${G'}^2$ can be extended to a 6-coloring of $G^2$.
Since, this is the lone remaining case for $\Delta=4$, this
yields the lemma. 

{\sc First case $k\in \{4,5\}$}:
For $k=4$ we first color $u_1$ and $u_3$ with a same color
that is unused at $v_1,v_4,x,y$. Then color the vertices $v_2,v_3,u_2$ 
in a greedy fashion in this order, using at most 6 colors.

For $k=5$ we first color $u_1$ and $v_4$ with a same color
that is unused at $v_1,v_5,x,y$. Then color $v_2$ and $u_4$ with a same
color that is unused at $v_1,v_5,x,y,u_1$ and $v_4$ (at most five colors used on these
six vertices). Finally color $v_3,u_2,u_3$ in a greedy fashion in this order, 
using at most 6 colors.

{\sc Second case $k\geq 6$:}
Note that in a coloring of $G^2$ in the previous case, $v_1,u_1,v_2,u_2,v_3$ 
all receive distinct colors since they
induce a clique in $G^2$. Assume these colors are $1,2,3,4,5$ 
in this order respectively.
Such a coloring can now be extended to a case of $k\geq 6$ 
by coloring $v_1,u_1,v_2,u_2,v_3,u_3,v_4,u_4,v_5,u_5,v_6$
using the colors $1,\ldots,6,1,\ldots 5$ in this order respectively. 
Since each case of $k\geq 6$ can be obtained by repeated use of such an 
extension, we have the lemma.
\end{proof}


\subsection{Cases with $\Delta=6$} 

We now delve into the case where $G$ is an outerplanar graph with $\Delta = 6$,
and we will show that 
$\chi(G^2) = \Delta + 1 = 7$ always holds here. By 
Lemma~\ref{lmm:biconn} we shall assume $G$ to be biconnected and hence,
by Lemma~\ref{lmm:dual-out}, its weak dual $T^*(G)$ to be a connected tree.
As in Section~\ref{sec:ind}, we will reduce our considerations to 
some key cases regarding the weak dual $T^*(G)$ of $G$. We will assume,
unless otherwise stated, that $G$ is a biconnected outerplanar plane graph
with $\Delta = 6$ with $T^*(G)$ as a connected tree. 

Assuming $\chi(G^2) > 7$, and that $G$ is minimal with this property
(that is, any other graph with $\Delta = 6$ has its square $7$-colorable),
we note that we can assume that we do not have configurations (A) and (B)
shown in Figure~\ref{fig:red1} nor configuration (C) in Figure~\ref{fig:red2},
(since in those cases we would have a $6$-vertex, contradicting
our assumption on $G$.). In other words, we may assume that
(i) each face $f$ with $f^*$ a leaf in $T^*(G)$ is bounded by three or four edges,
(ii) each 1-ss face $f$ is bounded by exactly three vertices, and
(iii) each 1-ss face $f$ has at most on sibling in $T^*(G)$.
We say that $G$ is \defn{{\em 3-restricted}} if it 
satisfies all these assumptions (i), (ii) and (iii).
In addition we have the following for the parent $f'$ of $f$ in $T^*(G)$.
\begin{lemma}
\label{lmm:parent-3-4}
If $G$ is 3-restricted, $f$ a 1-ss face and $|f'|\geq 5$,
then either $f$ or $f'$ contains a $6$-vertex.
\end{lemma}
\begin{proof}
Note that $f'$ (that is to say ${f'}^*$) is a leaf in the pruned tree 
$\pr(T^*(G))$. Let
$f'$ be bounded by the vertices $v_0,\ldots, v_{\alpha}$ with $\alpha\geq 4$
and where $v_0v_{\alpha}$ it the edge in $G$ that boarders the grandparent
$f''$ (that is to say, is dual to 
the edge in $T^*(G)$ incident to $f'$).
If $v_i$ and $v_{i+1}$ are on the boundary of either $f$ or its unique sibling
$g$ (in the case that $f$ has a sibling $g$) for some $i\in\{1,\ldots,\alpha - 2\}$, 
then the degree-two vertex on either $f$ or $g$ has degree at most five in $G^2$. 
Otherwise, all the $\alpha-2\geq 2$ edges 
$v_1v_2,\ldots,v_{\alpha-2}v_{\alpha-1}$ bound the infinite face of $G$,
in which case the $\alpha -3\geq 1$ vertices, $v_2,\ldots,v_{\alpha-2}$, 
all are degree-two vertices with at most six neighbors in $G^2$.
This completes the proof.
\end{proof}
Consider further the case for a 3-restricted $G$ where $f$ is a 1-ss 
face of $G$, the face $f'$ is bounded by four vertices 
$v_0,v_1,v_2,v_3$ and $f$ is bounded by $u,v_1,v_2$ in such a way
that the edges $v_0v_1$, $v_1v_2$ and $v_2v_3$ all bound the
infinite face of $G - u$. In this
case both $v_1$ and $v_2$ have degree three in $G$ and hence the degree-two vertex 
$u$ has four neighbors in $G^2$. To avoid any vertices of degree $\leq 6$ in $G^2$,
we can therefore assume that if, for any 1-ss face $f$
with $f'$ bounded by the four vertices $v_0,v_1,v_2,v_3$, then the edge
$v_1v_2$ must bound the infinite face of $G$.

To make further restrictions, assume that $G$ is a biconnected outerplanar 
graph that is induced by the cycle $C_n$ on the vertices $\{u_1,\ldots,u_n\}$
in clockwise order. If $d_G(u_3) = 4$ and $u_3$ is adjacent to both $u_1$ and $u_5$ in $G$, 
then for any coloring of the square $G^2$ the vertices $u_1,\ldots, u_5$ 
must all receive distinct
colors, say $1,\ldots, 5$ respectively, since $N_G[u_3] = \{u_1,\ldots,u_5\}$ 
induces a clique in $G^2$. Consider the outerplanar graph $G'$ obtained by 
first removing both the edges $u_3u_4$ and $u_3u_5$ and then connecting
a new vertex $u_3'$ to each of the vertices $u_3,u_4$ and $u_5$. In this way
$G$ becomes the contraction of $G'$, namely $G = {G'}/u_3u_3'$. Note that 
if $G$ has a maximum degree of $\Delta=6$, then so does $G'$. In addition, given
the mentioned coloring of $G^2$ where $u_i$ has color $i$ for $1\leq i\leq 5$, 
then we can obtain a coloring
of ${G'}^2$ by retaining the colors of $u_i$ from $G^2$ for all $i\not\in\{2,3,4\}$,
and then assigning colors $3,2,4,3$ to
the vertices $u_2,u_3,u_3',u_4$ respectively.
\begin{definition}
\label{def:min-crim}
A biconnected outerplanar graph $G$ with maximum degree $\Delta = 6$
and a minimum number of vertices satisfying $\chi(G^2) = 8$ is called 
a {\em minimal criminal}.
\end{definition}
Clearly each minimal criminal must be 3-restricted.
What our discussion preceding the above definition means, in particular,
is the following additional property of a potential minimal criminal $G$.
\begin{theorem}
\label{thm:restricting-f-f'}
If $G$ is a minimal criminal, then $G$ has no degree-two vertices 
with $\leq 6$ neighbors in $G^2$. Further, let $f$ be a 
1-ss face of $G$.
\begin{enumerate}
  \item If $f$ has no sibling, then $f'$ is bounded by four vertices. Further,
all the faces $f$, $f'$ and $f''$ have exactly on vertex in common on their boundaries. 
  \item If $f$ has one sibling $g^*$, then $f'$ is bounded by three vertices.
Hence, all the faces $f$, $g$ and $f'$ are bounded by exactly three vertices and edges.
\end{enumerate}
\end{theorem}
\begin{proof}
If $f$ has no sibling and $f'$ is bounded by three vertices and edges,
then the degree-two vertex $u$ bounding $f$ has $\leq 6$ neighbors in $G^2$.
Then by minimality of $|V(G)|$, the square of $G-u$ can be colored by at most seven
colors, and hence so can $G^2$, since there is at least on color left for $u$
among the seven colors available. If further $f'$ is bounded by four vertices
and the faces $f$, $f'$ and $f''$ have no common vertex on their boundaries,
then as previously noted, $u$ has exactly four neighbors in $G^2$
since both neighbors of $u$ have degree three in $G$. 

If $f$ has one sibling $g$, and $f'$ is bounded by four vertices 
$v_0,v_1,v_2,v_3$, then the edge $v_1v_2$ bounds neither $f$ nor
$g$ (since otherwise either $f$ or $g$ has a degree-two vertex on its
boundary with at most five neighbors in $G^2$) and therefore (assuming $f$ is
to the left of $g$ in the plane embedding of $G$) we have that
$v_0v_1$ bounds $f$, the edge $v_2v_3$ bounds $g$ and $v_1v_2$ bounds
the infinite face of $G$. Again, by minimality of $n$ we have that
the the square of the contracted graph $G/v_1v_2$ has a legitimate
7-coloring. By our above discussion preceding Definition~\ref{def:min-crim}, 
this coloring can then be extended to a 7-coloring of $G^2$, thereby 
contradicting the criminality of $G$. This complete the proof.
\end{proof}
What Theorem~\ref{thm:restricting-f-f'} implies, in particular, is that
in a 3-restricted minimal criminal $G$,
each configuration $C(f,f')$ of $f$ and its parent $f'$, where $f$ is
1-ss, is itself induced by a
5-cycle on the vertices $v_1,v_2,v_3,v_4,v_5$ in a clockwise order, 
and is of one of the following
three types (note that if $f$ has a sibling $g$, then it is unique and we assume
$g$ to be the right of $f$ in the planar embedding of $T^*(G)$ when viewed from $f'$):
\begin{description}
  \item[(a)] $C(f,f')$ is the 5-cycle on $\{v_1,v_2,v_3,v_4,v_5\}$
in which $v_3$ is connected to $v_1$. Here $f$ is bounded by $\{v_1,v_2,v_3\}$
and $f'$ is bounded by $\{v_1,v_3,v_4,v_5\}$. 
  \item[(b)] $C(f,f')$ is the 5-cycle on $\{v_1,v_2,v_3,v_4,v_5\}$
in which $v_3$ is connected to $v_5$. Here $f$ is bounded by $\{v_3,v_4,v_5\}$
and $f'$ is bounded by $\{v_1,v_2,v_3,v_5\}$. 
  \item[(c)] $C(f,f')$ is the 5-cycle on $\{v_1,v_2,v_3,v_4,v_5\}$
in which $v_3$ is connected to both $v_1$ and $v_5$. Here $f$ is bounded by $\{v_1,v_2,v_3\}$,
the face $g$ is bounded by $\{v_3,v_4,v_5\}$ and $f'$ is bounded by $\{v_1,v_3,v_5\}$. 
\end{description}
Here, for all the three types of configurations, it is assumed that 
the edge $v_1v_5$ bounds the face $f''$ as well as $f'$.

{\sc Remarks:}
Note that as plane configurations (a) and (b) are mirror images of each other.
Also, note that in all configurations, all the edges $v_iv_{i+1}$ where $1\leq i\leq 4$
of the 5-cycle that induces $C(f,f')$, except one edge $v_1v_5$,
bound the infinite face of $G$.

\vspace{3 mm}

{\bf Convention:} 
Let $G$ be a 3-restricted biconnected outerplanar graph with $\Delta = 6$
such that for each 1-ss face $f$ of $G$ the 
configuration of $f$ and its parent $f'$ is of type (a), (b) or (c) from
above. Call such a $G$ {\em fully restricted}.
Hence, a minimal criminal is always fully restricted. 

Let $G$ be a biconnected outerplanar graph induced by the cycle $C_n$ on
the vertices $\{u_1,\ldots,u_n\}$ in clockwise order. Assume that 
$d_G(u_4) = 6$ and that $u_4$ is adjacent to $u_1, u_2, u_3, u_5, u_6$ and $u_7$. 
From $G$ we construct four other outerplanar graphs $\tilde{G}$, $G'$, $G''$ and $G'''$ in the
following way:
\begin{enumerate}
  \item Let $\tilde{G}$ be obtained by replacing the edge $u_1u_2$
by the 2-path $(u_1,x,u_2)$.
  \item Let $G'$ be obtained from $G$ by replacing the edges $u_1u_2$
and $u_6u_7$ by the 2-paths $(u_1,x,u_2)$ and $(u_6,y,u_7)$
respectively.
  \item Let $G''$ be obtained from $G$ by replacing the edge $u_1u_2$
by the 2-path $(u_1,x,u_2)$ and connecting the additional vertex $y$
to both $u_6$ and $u_7$.
  \item Let $G'''$ be obtained from $G$ by connecting the additional
vertex $x$ to both $u_1$ and $u_2$ and the additional vertex $y$ to $u_6$ and $u_7$.
\end{enumerate}
Note that $G$ is a contraction of each of the graphs $\tilde{G}$, $G'$, $G''$ and $G'''$,
namely 
\[
\tilde{G}/u_1x = ({G'}/u_1x)/u_6y = ({G''}/u_1x)/u_6y 
                    = ({G'''}/u_1x)/u_6y = G.
\]
\begin{lemma}
\label{lmm:interm-col-ext}
A 7-coloring of $G^2$ can be extended to a 7-coloring of ${\tilde{G}}^2$.
\end{lemma}
\begin{proof}
Since $N_G[u_4] = \{u_1, u_2, u_3, u_4, u_5, u_6, u_7\}$ induces
a clique in $G^2$, we may assume $u_i$ to have color $i$ for
$1\leq i\leq 7$ in the 7-coloring of $G^2$. To obtain a 7-coloring of ${\tilde{G}}^2$
we retain the colors of $u_i$ from $G^2$ for all $i\not\in\{2,3,5\}$ and 
then assign colors $2,3,5,2$ to vertices $x,u_2,u_3,u_5$ respectively
(note that we do not need to know the colors
of all the neighbors of neither $u_1$ nor $u_7$ in the given 7-coloring of $G^2$).
\end{proof}
{\bf Convention:} 
Let $G$ be a biconnected outerplanar graph $G$ with $\Delta = 6$.
Let $u$ be a degree-two vertex of $G$ that bounds a leaf-face $f$ 
of $T^*(G)$. If $(G - u)^2$ is provided with a 7-coloring, 
call $u$ {\em c-simplicial} if all the neighbors of $u$ in $G^2$
have collectively $\leq 6$ colors.

Note that $\tilde{G}$ from above cannot be a minimal criminal, since $u_5$
is c-simplicial; if we have a 7-coloring of $(\tilde{G} - u_5)^2$ then we
can extend it to a 7-coloring of ${\tilde{G}}^2$.

Our next theorem will provide our main tool for this section.
\begin{theorem}
\label{thm:main-re-coloring}
If $G$ and the constructed graphs $G'$, $G''$ and $G'''$ are as defined above,
then none of the graphs $G'$, $G''$ or $G'''$ are minimal criminals.
\end{theorem}
\begin{proof}
If $G'$ is a minimal criminal, then by definition $G^2$ has a legitimate 7-coloring. 
Again, we may assume $u_i$ to have color $i$ for $1\leq i\leq 7$. 
By retaining the colors of $u_i$ from $G^2$ for $i\not\in\{2,3,5,6\}$ 
and then assigning colors $2,3,6,2,5,6$ to the vertices $x,u_2,u_3,u_5,u_6,y$
respectively, we obtain a legitimate 7-coloring of ${G'}^2$. 
Hence, $G'$ cannot be a minimal criminal.

If $G''$ is a minimal criminal, then $G^2$ has a legitimate 7-coloring.
We can assume $u_i$ to have color $i$ for $1\leq i\leq 7$. 
By Lemma~\ref{lmm:interm-col-ext} we obtain a 7 coloring of ${\tilde{G}}^2$,
as given in its proof. If $y$ is c-simplicial (w.r.t.~this mentioned 7-coloring of $G^2$)
then we can obtain a 7-coloring of ${G''}^2$. Therefore $y$ cannot be c-simplicial
in this case. This means the neighbors of $u_7$ among $V(G)\setminus\{u_2,\ldots,u_6\}$
have the colors $1,3$ and $5$ precisely, since $d_G(u_7) = 5$ and $d_{G''}(u_7) = 6$.
In this case assign the colors $2,5,6,3,2,6$ to the vertices $x,u_2,u_3,u_5,u_6,y$. 
This is a legitimate 7-coloring of ${G''}^2$ 
and hence $G''$ cannot be a minimal criminal.

If $G'''$ is a minimal criminal, then then $G^2 = (G'' -  xy)^2$ has a legitimate 7-coloring.
We can assume $u_i$ to have color $i$ for $1\leq i\leq 7$. If both $x$ and $y$
are c-simplicial, then we can extend the given coloring of $G^2$ to that of ${G'''}^2$,
since $x$ and $y$ are of distance 3 or more from each other in $G'''$. If neither
$x$ nor $y$ are c-simplicial, then we must have that the neighbors of $u_1$ among 
$V(G)\setminus\{u_2,\ldots,u_6\}$ have the colors $5,6$ and $7$ precisely,
and the neighbors of $u_7$ among $V(G)\setminus\{u_2,\ldots,u_6\}$
have the colors $1,2$ and $3$ precisely. In this case we assign the colors
colors $2,3,6,2,5,6$ to the vertices $x,u_2,u_3,u_5,u_6,y$
respectively (as in the case with $G'$) and obtain a legitimate 7-coloring of ${G'''}^2$. 
We consider lastly the case where one of $x$ and $y$ is c-simplicial and the other is not. 
By symmetry, it suffices to consider the case where $x$ is c-simplicial and $y$ is not.
The fact that $x$ is c-simplicial means that it can be assigned a color that must be from 
$\{5,6,7\}$ and thereby obtain a 7-coloring of $G''' - y$. Since $y$ is not
c-simplicial means that the neighbors of $u_7$ among $V(G)\setminus\{u_2,\ldots,u_6\}$
have the colors $1,2$ and $3$ precisely. We now consider the following three
cases:

{\sc $x$ has color $5$:} In this case assign the colors
$2,5,3,2,6,5$ to the vertices $x,u_2,u_3,u_5,u_6,y$, thereby obtaining a legitimate
7-coloring of ${G'''}^2$.

{\sc $x$ has color $6$:} In this case assign the colors
$2,6,3,2,5,6$ to the vertices $x,u_2,u_3,u_5,u_6,y$, thereby obtaining a legitimate
7-coloring of ${G'''}^2$.

{\sc $x$ has color $7$:} Here $u_1$ and $u_7$ cannot be connected since
both $x$ and $u_7$ have color $7$. In this case assign the colors
$7,2,5,3,6,5$ to the vertices $x,u_2,u_3,u_5,u_6,y$, thereby obtaining a legitimate
7-coloring of ${G'''}^2$.

This shows that $G'''$ cannot be a minimal criminal. This completes our proof.
\end{proof}
{\sc Remark:} To test the legitimacy of the extended colorings we note
first of all that the vertices $u_1$, $u_4$ and $u_7$ always keep their
color from the one provided by $G^2$. In addition, the colors of the
neighbors of $u_1$ among $\{u_1,\ldots,u_7\}$ are the same, unless $x$
is not c-simplicial, which gives concrete information about the colors
of the other neighbors of $u_1$. Similarly the colors of the neighbors of $u_7$ 
among $\{u_1,\ldots,u_7\}$ are the same, unless (as for $x$) 
$y$ is not c-simplicial, which again gives concrete information about the
colors of the remaining neighbors of $u_7$.

\vspace{3 mm}

We are now ready for the proof of the following main result of this section.
\begin{theorem}
\label{thm:main-D=6}
There is no minimal criminal; every biconnected outerplanar graph $G$ with
$\Delta = 6$ has $\chi(G^2) = \Delta + 1 = 7$.
\end{theorem}
\begin{proof}
We will show that a minimal criminal must have the form of one of the
graphs $G'$, $G''$ or $G'''$, thereby obtaining a contradiction by
Theorem~\ref{thm:main-re-coloring}.

Assume $G$ is a minimal criminal, which must therefore be fully restricted.
Since $\Delta = 6$, each 1-ss face $f$ of $G$ has a 
parent $f'$ and a grandparent $f''$ in $T^*(G)$. Hence, $f''$ (that is ${f''}^*$)
is a leaf in $\pr^2(T^*(G))$ (or a single vertex). Since $G$ is fully 
restricted, the configuration $C(f,f')$ in $G$ is induced by a 5-cycle
on vertices $\{v_1,v_2,v_3,v_4,v_5\}$ in clockwise order and is of type (a), (b) or (c) 
mentioned earlier. Assume $f''$ is bounded by $u_1,\ldots,u_m$
where $m\geq 3$. If $f''$ (that is ${f''}^*$) is not a single vertex but a leaf in 
$\pr^2(T^*(G))$, then let the edge $u_mu_1$ of $G$ be the dual edge
of the unique edge with ${f''}^*$ as and endvertex in $\pr^2(T^*(G))$.
In any case (whether ${f''}^*$ is a leaf or a single vertex in $\pr^2(T^*(G))$) 
at least one of the edges $u_1u_2, \ldots,u_{m-1}u_m$ must be
identified with an edge $v_1v_5$ of a configuration $C(f,f')$ 
of type (a), (b) or (c).

If there is an edge $u_iu_{i+1}$ bounding $f''$ and a configuration
$C(f,f')$ in such a way that the edge $v_1v_5$ is identified with
$u_iu_{i+1}$ (i.e.~$v_1 = u_i$ and $v_5 = u_{i+1}$) and such
that either $d_G(v_1) \leq 5$ or $d_G(v_5) \leq 5$, then either the degree-two
vertices $v_1$ or $v_5$ has $\leq 6$ neighbors in $G^2$ respectively.
this means that a 7-coloring of either $G/v_1v_2$ or $G/v_4v_5$ can 
be used to extend to a 7-coloring of $G$. Therefore $G$ cannot be a minimal
criminal in this case. 

We note that in order for $d_G(v_1)=d_G(v_5)=6$ for all configurations
$C(f,f')$, then {\em every} edge $u_iu_{i+1}$ for $1\leq i\leq m$ must
be identified with an edge $v_1v_5$ of a configuration $C(f,f')$.
That is to say, none of the edges $u_iu_{i+1}$ can bound the infinite
face of $G$. Since $m\geq 3$ we must, in particular, have that 
the edges $u_1u_2$ and $u_2u_3$ must be identified
with edges $v_1v_5$ of configurations $C(f,f')$ each of type
(a), (b) or (c). If $d_G(u_2)\leq 5$, then (as mentioned in previous
paragraph) both configurations, to the left and right of $u_2$ in the
plane embedding of $G$, contain a degree-two vertex with $\leq 6$ 
neighbors in $G^2$. In order for $d_G(u_2)=6$  then the $C(f,f')$
configuration to the left of $u_2$ in $G$, must be of type (b) or (c)
and the configuration $C(f,f')$ to the right of $u_2$ must be of
type (a) or (c). We now finally discuss these cases separately.
Here ``{\sc Case {\rm (x,y)}}'' will mean that configuration $C(f,f')$ of type
(x) is to the left of $u_2$ and configuration $C(f,f')$ of type (y)
is to the right of $u_2$.

{\sc Case {\rm (b,a)}:} Here $G$ is of type $G'$ as stated in 
Theorem~\ref{thm:main-re-coloring} (with $u_2$ in the role of $u_4$ mentioned
there), and therefore cannot be a minimal criminal.

{\sc Case {\rm (b,c)}:} Here $G$ is of type $G''$ as stated in 
Theorem~\ref{thm:main-re-coloring}, and therefore cannot be a minimal criminal.

{\sc Case {\rm (c,a)}:} Here $G$ is a mirror image of a type $G''$ (previous case)
as stated in Theorem~\ref{thm:main-re-coloring}, 
and therefore cannot be a minimal criminal.

{\sc Case {\rm (c,c)}:} Here $G$ is of type $G'''$ as stated in 
Theorem~\ref{thm:main-re-coloring}, and therefore cannot be a minimal criminal.

This concludes the proof, that there is no minimal criminal, and hence
the square of each biconnected outerplanar graph with $\Delta = 6$ is 7-colorable.
\end{proof}
By Theorem~\ref{thm:main-D=6} and Lemma~\ref{lmm:biconn} we have the following corollary.
\begin{corollary}
\label{cor:main-D=6}
For every outerplanar graph $G$ with $\Delta = 6$ we have $\chi(G^2) = 7$.
\end{corollary}


\subsection{Greedy is not exact}

We have observed that {\Greedy} yields an optimal for coloring squares of
outerplanar graphs whenever $\Delta \ge 7$ or $\Delta=2$.
On many of the examples that we have constructed it also gives optimal
colorings. It is therefore a natural question to ask whether it always
obtains an optimal coloring. If not, one may ask for the case of
chordal graphs, where we have seen that {\Greedy} is also optimal for $\Delta=3$ and $\Delta=4$.
Further evidence may be gathered by observing that it yields optimal
colorings of $F_4$, $F_5$ and $RL_n$, since their squares are chordal.

We answer these questions in the negative by showing that the chordal
graph $F_6$ is a counterexample.

\begin{theorem}
{\Greedy} does not always output an optimal coloring of $F^2_{6}$.
\end{theorem}

\begin{proof}
Suppose that the vertices of (\ref{fig:subfig:f6}) in Figure~\ref{fig:f456} 
are ordered 
so that first come the white vertices in the order shown on the figure
(either from left-to-right or in a circular order). Then, Greedy will
first color the white vertices with the first two colors. Now, the
remaining six vertices must receive different colors, and none of them
can use the first two colors, resulting in an 8-coloring.
\end{proof}


\subsection*{Acknowledgments}  

The authors are grateful to Steve Hedetniemi
for his interest in this problem and for suggesting
the writing of this article.



\begin{thebibliography}{10}


\bibitem{GeirMagnus}
\newblock G.~Agnarsson, M.~M.~Halld\'{o}rsson.
\newblock Coloring Powers of Planar Graphs.
\newblock \emph{SIAM Journal of Discrete Mathematics}, \textbf{16}, No.~4, 651 -- 662, (2003).

\bibitem{GeirMagnus-SODA2}
\newblock G.~Agnarsson, M.~M.~Halld\'{o}rsson.
\newblock On Colorings of Squares of Outerplanar Graphs.
\newblock {\em Proceedings of the Fifteenth Annual ACM-SIAM Symposium On Discrete
               Algorithms (SODA 2004), New Orleans}, 237 -- 246, (2004).


\bibitem{AppelHaken}
\newblock K.~Appel, W.~Haken. 
\newblock Every Planar Graph is Four Colorable.
\newblock Contemporary Mathematics 98, AMS, Providence, (1989). 



\bibitem{BBGH01}
\newblock O.~Borodin, H.~J.~Broersmo, A.~Glebov, and J.~van~den~Heuvel.
\newblock Stars and bunches in planar graphs. Part II: General planar graphs and colorings. 
\newblock {\em CDAM Research Report Series} 2002-05, (2002).

\bibitem{CalaPet}
\newblock T. Calamoneri and R. Petreschi.
\newblock $L(h,1)$-labeling subclasses of planar graphs.
\newblock {\em Journal of Parallel and Distributed Computing}, \textbf{64}(3):414--426, (2004).


\bibitem{Diestel}
\newblock R.~Diestel.
\newblock Graph Theory.
\newblock {\em Graduate Texts in Mathematics}, GTM -- 173, Springer Verlag, (1997).


\bibitem{mmh:greed}
\newblock M. M. Halld{\'o}rsson and J. Radhakrishnan.
\newblock Greed is good: Approximating independent sets in sparse and bounded-degree graphs.
\newblock {\em Algorithmica}, \textbf{18}: 145 -- 163, (1997).

\bibitem{JensenToft}
\newblock T.~R. Jensen and B.~Toft.
\newblock Graph Coloring Problems.
\newblock {\em Wiley Interscience}, (1995).
\newblock {\tt http://www.imada.sdu.dk/Research/Graphcol/}.


\bibitem{KMR:dialm}
\newblock S.~O. Krumke, M.~V. Marathe, and S.~S. Ravi.
\newblock Approximation algorithms for channel assignment in radio networks.
\newblock In {\em Dial M for Mobility, 2nd International Workshop on Discrete
  Algorithms and Methods for Mobile Computing and Communications}, Dallas,
  Texas, (1998).

\bibitem{LaskarShier}
\newblock R.~Laskar and D.~Shier.
\newblock On Powers and Centers of Chordal Graphs.
\newblock {\em Discrete Applied Mathematics}, \textbf{6}: 139 -- 147, (1983).


\bibitem{MS01}
\newblock M.~Molloy and M.~R.~Salavatipour.
\newblock A bound on the chromatic number of the square of a planar graph.
\newblock {\em J.~Combinatorial Theory (Series B)}, \textbf{94}(2):189--213, (2005). 




\bibitem{Wegner}
\newblock G.~Wegner.
\newblock Graphs with given diameter and a coloring problem.
\newblock {\em Technical report}, University of Dortmund, (1977).

\bibitem{West}
\newblock D.~B.~West.
\newblock Introduction to Graph Theory.
\newblock {\em Prentice-Hall Inc.}, 2nd ed., (2001).

\bibitem{ZKN00}
\newblock X.~Zhou, Y.~Kanari, and T.~Nishizeki.
\newblock Generalized vertex-colorings of partial $k$-trees.
\newblock {\em IEICE Trans. Fundamentals}, E83-A: 671 -- 678, (2000).

\end{thebibliography}
\end{document}